\documentclass{article}[12pt]

\usepackage{amssymb, amsmath, amsthm,commath,stackrel,enumerate,mathtools,bbold,xcolor,thm-restate,geometry}

\author{Jeroen Winkel}
\date{\today}
\title{Geometric property (T) for non-discrete spaces}

\theoremstyle{definition}
\newtheorem{example}{Example}[section]
\theoremstyle{plain}

\newtheorem{lemma}[example]{Lemma}
\newtheorem{proposition}[example]{Proposition}

\newtheorem{corollary}[example]{Corollary}

\newtheorem{question}[example]{Question}
\theoremstyle{remark}
\newtheorem{remark}[example]{Remark}
\theoremstyle{definition}
\newtheorem{definition}[example]{Definition}

\renewcommand{\phi}{\varphi}
\renewcommand{\H}{\mathcal H}
\renewcommand{\epsilon}{\varepsilon}
\DeclareMathOperator{\cs}{cs}

\DeclareMathOperator{\im}{im}
\DeclareMathOperator{\tensor}{\otimes}
\DeclareMathOperator{\Warp}{Warp}

\begin{document}
\maketitle

\begin{abstract}
    Geometric property (T) was defined by Willett and Yu, first for sequences of graphs and later for more general discrete spaces.
    Increasing sequences of graphs with geometric property (T) are expanders, and they are examples of coarse spaces for which the maximal coarse Baum-Connes assembly map fails to be surjective.
    Here, we give a broader definition of bounded geometry for coarse spaces, which includes non-discrete spaces.
    We define a generalisation of geometric property (T) for this class of spaces and show that it is a coarse invariant.
    Additionally, we characterise it in terms of spectral properties of Laplacians.
    We investigate geometric property (T) for manifolds and warped systems.
\end{abstract}

\section{Introduction}
In their paper \cite{WilYu12}, Willett and Yu introduce \emph{geometric property (T)} for sequences of graphs, in order to provide examples of coarse spaces for which the maximal coarse Baum-Connes assembly map fails to be surjective.
The concept is studied further by the same authors in \cite{WilYu14}.
There, geometric property (T) is defined for coarse spaces that are discrete with bounded geometry, monogenic, and countable.
Such a space has geometric property (T) if each unitary representation of the translation algebra $\mathbb C_u[X]$ (we call it $\mathbb C_{\cs}[X]$, see Definition \ref{def-algebra-controlled-support}) that has almost invariant vectors, has a non-zero invariant vector.
For increasing sequences of finite graphs, this is a strictly stronger property than being an expander.
In fact, in this case it is equivalent to the graph Laplacian having spectral gap in the maximal completion of $\mathbb C_u[X]$ (see \cite[Proposition 5.2]{WilYu14}).
A sequence of finite graphs with geometric property (T) can never have large girth (see \cite[Corollary 7.5]{WilYu12}).

The aim of this article is to generalise geometric property (T) to spaces satisfying a broader notion of \emph{bounded geometry}, see Definition \ref{def-geometric-T}.
This broader notion of bounded geometry is explained in Section \ref{sec-bounded-geometry}.
Spaces with bounded geometry include sequences of graphs with bounded degree, manifolds whose Ricci curvature and injectivity radius are uniformly bounded from below, and warped systems coming from an action of a finitely generated group on a compact manifold.
We do not assume monogenicity of the spaces.
However, we do assume that the spaces are covered by a countable number of ``bounded" sets.

In \cite{WilYu14}, \emph{translation operators} are used extensively to study geometric property (T) for discrete spaces.
For non-discrete spaces, these operators do not behave as nicely.
Therefore we introduce the concept of ``operators in blocks", which are operators with a convenient decomposition (see Definition \ref{def-operators-in-blocks}).

Any coarse space with bounded geometry can be equipped with a measure $\mu$ that is \emph{uniformly bounded} and for which a \emph{gordo} set exists (see Definition \ref{def-uniformly-bounded-gordo} and Proposition \ref{prop-bg-coarse}).
We will firstly define geometric property (T) for a coarse space $X$ equipped with such a measure $\mu$.
However, Theorem \ref{thm-coarse-invariance} will show that it is in fact independent of the chosen measure $\mu$.

We will consider the Hilbert space $L^2(X,\mu)$.
An operator $T$ in $B(L^2(X,\mu))$ has \emph{controlled support} if the support of a function is not changed much by $T$ (see Definition \ref{def-controlled-support}).
Let $\mathbb C_{\cs}[X]\subseteq B(L^2(X,\mu))$ be the algebra of operators with controlled support.
A representation of this algebra is a unital $*$-homomorphism $\rho\colon \mathbb C_{\cs}[X]\to B(\H)$ for some Hilbert space $\H$.
For such a representation, the subspace of \emph{constant vectors} $\H_c$ is defined (see Definition \ref{def-constant-vector}).
The coarse space $X$ has geometric property (T) if every unit vector $v\in \H_c^\perp$ is ``moved enough" by a suitable operator in $\mathbb C_{\cs}[X]$ (see Definition \ref{def-geometric-T}).

An important notion in coarse geometry is that of \emph{coarse invariance} (see Section \ref{sec-coarse-spaces}).
In \cite{WilYu14}, it was shown that geometric property (T) is a coarse invariant (see \cite[Theorem 4.1]{WilYu14}).
In this paper we show that our generalisation of geometric property (T) is still a coarse invariant.
This gives evidence that we have the ``right" definition of property (T).
\begin{restatable*}{theorem}{coarseinvariance}\label{thm-coarse-invariance}
Suppose $(X,\mu)$ and $(X',\mu')$ are coarsely equivalent spaces with bounded geometry, equipped with uniformly bounded measures for which gordo sets exist.
Then $(X,\mu)$ has property (T) if and only if $(X',\mu')$ has property (T).
In particular, whether $X$ has property (T) is independent of the chosen uniformly bounded measure $\mu$ for which a gordo set exists.
\end{restatable*}

Analogous to \cite{WilYu14}, we establish a complete characterisation of geometric property (T)  in terms of non-amenability for connected coarse spaces of bounded geometry.
In fact, for unbounded connected spaces, geometric property (T) is the complete converse of amenability for coarse spaces.
This generalizes \cite[Corollary 6.1]{WilYu14}.
The notion of amenability for coarse spaces with bounded geometry is introduced in Lemma \ref{prop-amenable-coarse-space}.
\begin{restatable*}{theorem}{amenability}\label{thm-connected-space-amenable-or-(T)}
Let $X$ be a connected coarse space and let $\mu$ be a uniformly bounded measure for which a gordo set exists.
If $X$ is bounded, then it is amenable and has geometric property (T).
If $X$ is unbounded, it has geometric property (T) if and only if it is not amenable.
\end{restatable*}

Towards the end of this article, we investigate geometric property (T) for Riemannian manifolds.
Let $M$ be a Riemannian manifold or a countable disjoint union of Riemannian manifolds with the same dimension.
The Riemannian metric defines a coarse structure on $M$.
If the Ricci curvature and the injectivity radius of $M$ are uniformly bounded from below, then $M$ has bounded geometry (see Section \ref{sec-manifolds}).
In this case the manifold Laplacian $\Delta_M$ is an element of the algebra $C^*_{\max}(X)$, and the manifold $M$ has geometric property (T) if and only  if $\Delta_M$ has spectral gap in this algebra.

\begin{restatable*}{theorem}{manifolds}\label{thm-manifold-(T)-iff-Laplacian-spectral-gap}
Let $(M,g)$ be a Riemannian manifold of bounded geometry equipped with the coarse structure coming from the geodesic metric.
Then $M$ has geometric property (T) if and only if the Laplacian $\Delta_M$ has spectral gap in $C^*_{\max}(M)$.
\end{restatable*}

In the last section we discuss warped systems.
A group $\Gamma$ acting by homeomorphisms on a Riemannian manifold $M$ gives rise to a coarse spaces called a \emph{warped system}, see Section \ref{sec-warped-systems}.
In Section \ref{sec-warped-systems} we will show that warped systems have bounded geometry if the manifold $M$ is compact.
One might expect that the warped system has geometric property (T) if and only if the group has property (T) and the action is ergodic.
This remains an open question, though some partial results are given.

\section{Preliminaries on coarse spaces}\label{sec-coarse-spaces}
A coarse space is a space with a large-scale geometric structure.
They were first described in \cite{HigPedRoe95} and further developed in \cite{Roe03}.
Let us first recall the definition from \cite[Definition 2.3]{Roe03}.

\begin{definition}
Let $X$ be a set.
A \emph{coarse structure} on $X$ consists of a collection $\mathcal E$ of subsets of $X\times X$, whose elements are called the \emph{controlled sets}, such that: 
\begin{enumerate}[(i)]
    \item The diagonal $\{(x,x)\mid x\in X\}$ is an element of $\mathcal E$.
    \item For every $E\in\mathcal E$ and $F\subseteq E$, we have $F\in\mathcal E$.
    \item For every $E\in\mathcal E$, the \emph{inverse} $E^{-1}=\{(x,y)\in X\times X\mid (y,x)\in E\}$ is an element of $\mathcal E$.
    \item For every $E,F\in\mathcal E$, we have $E\cup F\in\mathcal E$.
    \item For every $E,F\in\mathcal E$, the \emph{composition} 
    \[E\circ F=\{(x,z)\in X\times X\mid \text{there is }y\in X\text{ such that }(x,y)\in E\text{ and }(y,z)\in F\}\]
    is an element of $\mathcal E$.
\end{enumerate}
A set $X$ equipped with a coarse structure is called a \emph{coarse space}.
\end{definition}

\begin{example}\label{ex-metric-coarse-space}
Every metric space $(X,d)$ can be given a coarse structure:  let
\[\mathcal E=\{E\subseteq X\times X\mid \sup_{(x,y)\in E}d(x,y)<\infty\}.\]
For any $R>0$, we denote $E_R=\{(x,y)\in X\times X\mid d(x,y)\leq R\}$.
Then the controlled sets are exactly the subsets of $X\times X$ that are contained in some $E_R$.
\end{example}

\begin{definition}
Let $X$ and $Y$ be coarse spaces.
A map $f\colon X\to Y$ is a \emph{coarse equivalence} if there is a map $g\colon Y\to X$ such that for any controlled set $E\subseteq X\times X$ the set $(f\times f)(E)$ is controlled, for any controlled set $F\subseteq Y\times Y$ the set $(g\times g)(F)$ is controlled, the set $\{(x,gf(x))\mid x\in X\}$ is controlled and the set $\{(y,fg(y))\mid y\in Y\}$ is controlled.
The spaces $X$ and $Y$ are called \emph{coarsely equivalent} if such an $f$ exists.
\end{definition}
This definition is easily seen to be equivalent to Definition 2.21 in \cite{Roe03}.

An important example of the above is when $Y$ is coarsely dense in $X$.
\begin{definition}
Let $X$ be a coarse space and $Y\subseteq X$.
The subset $Y$ is called \emph{coarsely dense in $X$} if there is a controlled set $E\subseteq X\times X$ such that for all $x\in X$ there is $y\in Y$ with $(x,y)\in E$.
\end{definition}

If $Y$ is coarsely dense in $X$, it is easy to show that $X$ is coarsely equivalent to $Y$ with the induced coarse structure.

In this paper we will use the following properties and notation:
\begin{definition}\label{def-defs-coarse-set}
Let $X$ be a coarse space and $E\subseteq X\times X$ a controlled set.
\begin{enumerate}[(i)]
    \item The controlled set $E$ is called \emph{symmetric} if $E^{-1}=E$.
    \item For any $x\in X$, we write $E_x=\{y\in X\mid (y,x)\in E\}$.
    \item For any $U\subseteq X$, we write $E_U=\bigcup_{x\in U}E_x$.
    \item A subset $U\subseteq X$ is called \emph{bounded} if $U\times U$ is a controlled set.
    \item A subset $U\subseteq X$ is called \emph{$E$-bounded} if $U\times U\subseteq E$.
    \item We write $E^{\circ n}$ for the $n$-fold composition $E\circ E\circ \cdots \circ E$.
    \item The set $E$ \emph{generates} the coarse structure if for any controlled set $F$, there is an integer $n$ such that $F\subseteq E^{\circ n}$.
\end{enumerate}
\end{definition}

For later use we give some basic properties using these definitions.
\begin{lemma}\label{lem-basic-properties-coarse-spaces}
Let $X$ be a coarse space.
\begin{enumerate}[(i)]
    \item Each controlled set is contained in a symmetric controlled set.
    \item A subset $U\subseteq X$ is bounded if and only if it is $E$-bounded for some controlled set $E$.
    \item If $E$ is a controlled set and $x\in X$ then $E_x$ is $E\circ E^{-1}$-bounded (though it is not necessarily $E$-bounded).
    \item If $U$ is an $E$-bounded set, it is contained in $E_x$ for any $x\in U$.
\end{enumerate}
\end{lemma}
\begin{proof}
These follow directly from the definitions.
\end{proof}
Most of the time, we will only consider symmetric controlled sets.

\section{Bounded geometry}\label{sec-bounded-geometry}

Consider a coarse space $X$.
In \cite{WilYu14} it is defined that $X$ is a \emph{discrete space of bounded geometry} if for every controlled set $E\subseteq X\times X$, there is an integer $N$ such that $\#E_x\leq N$ for all $x\in X$.
We will define geometric property (T) for a wider class of spaces.
We still need a concept of bounded geometry, which we define below.
\begin{definition}\label{def-covering}
Let $X$ be a coarse space.
A controlled set $F\subseteq X\times X$ is called \emph{covering} if it is symmetric and for every controlled $E$ there is an $N$ such that each $E_x$ can be covered by at most $N$ sets that are $F$-bounded.
\end{definition}

\begin{lemma}\label{lem-different_def_covering}
Let $X$ be a coarse space and $F\subseteq X\times X$ a symmetric controlled set.
The following are equivalent:
\begin{enumerate}[(i)]
    \item The controlled set $F$ is covering.
    \item For each controlled set $E$ there is a constant $N$ such that each $E$-bounded set $U$ can be covered by at most $N$ sets that are $F$-bounded.
\end{enumerate}
\end{lemma}
\begin{proof}
Suppose $F$ is covering and let $E$ be a controlled set.
Let $N$ be such that each $E_x$ can be covered by at most $N$ sets that are $F$-bounded.
Let $U$ be a non-empty $E$-bounded set.
For any $x\in U$ we have $U\subseteq E_x$, so $U$ is also covered by at most $N$ sets that are $F$-bounded.

Conversely, suppose that $(ii)$ is true and let $E$ be a controlled set.
There is a constant $N$ such that each $E\circ E^{-1}$-bounded set can be covered by at most $N$ sets that are $F$-bounded.
Then each $E_x$ can also be covered by at most $N$ sets that are $F$-bounded.
\end{proof}

\begin{definition}\label{def-bounded-geometry}
A space $X$ has bounded geometry if there is a controlled covering set $F\subseteq X\times X$.
\end{definition}

The definition above is easily seen to be equivalent to the definition given by Roe in \cite[Definition 3.8]{Roe03}.

\begin{example}
Consider a metric space $X$ with the corresponding coarse structure.
If there exists a controlled covering set, we can enlarge it to $E_R$ for some $R>0$.
Therefore, $X$ has bounded geometry if and only if $E_R$ is covering for some $R>0$.
\end{example}
In fact, for many examples of metric spaces, $E_R$ will be covering for any $R>0$.

\begin{example}\label{ex-discrete-space-bounded-geometry}
A discrete space of bounded geometry certainly has bounded geometry in this definition, as we can take $F$ to be the diagonal.
\end{example}

A different way to characterise bounded geometry is by the existence of a measure on $X$ satisfying certain properties.
We will study measures on any sigma algebra on the set $X$, with as only condition, that singletons must be measurable.
A controlled set $E\subseteq X$ acts on the set of subsets of $X$, by sending $U$ to $E_U$.
Accordingly, given a measure on $X$, we say that a controlled set $E$ is \emph{measurable} if $E_U$ is measurable for all measurable $U\subseteq X$.
In particular, $E_x$ is then measurable.
Note that $(E\circ F)_U =E_{F_U}$, therefore the composition of measurable controlled sets is again measurable.

\begin{definition}\label{def-uniformly-bounded-gordo}
Let $(X,\mathcal E)$ be a coarse space, and let $\mu$ be a measure on $X$, for any sigma algebra that contains all singletons.
\begin{enumerate}[(i)]
\item The measure is called \emph{$\mathcal E$-uniformly bounded} or simply \emph{uniformly bounded} if for each controlled set $E$ there is a constant $C>0$ such that each $E$-bounded measurable set $U$ has measure at most $C$.
\item A symmetric controlled set $E$ is called \emph{$\mu$-gordo} or simply \emph{gordo} if it is measurable and $\mu(E_x)$ is bounded away from zero independently of $x$, for all $x\in X$.
\end{enumerate}
\end{definition}

\begin{proposition}\label{prop-bg-coarse}
Let $X$ be a coarse space.
The following are equivalent:
\begin{enumerate}[(i)]
    \item The space $X$ has bounded geometry.
    \item There is a uniformly bounded measure $\mu$ on $X$ and a gordo set $E$ for $\mu$.
    \item There is a coarsely dense subspace $Y\subseteq X$ that is a discrete space of bounded geometry.
\end{enumerate}
\end{proposition}
\begin{proof}
Suppose $X$ has bounded geometry and $F\subseteq X\times X$ is a controlled covering set.
By Zorn's Lemma there is a maximal subset $Y\subseteq X$ satisfying $(y,y')\not\in F$ for all $y\neq y'\in Y$.
By maximality we know that for each $x\in X$ there exists $y\in Y$ such that $(x,y)\in F$, so $Y$ is coarsely dense in $X$.
To show that $Y$ is a discrete space of bounded geometry, consider a controlled set $E\subseteq X\times X$.
There exists $N$ such that every $E_y$ can be covered by $F$-bounded sets $U_1,\ldots,U_N$.
Since each $F$-bounded set can contain at most 1 element of $Y$, we see that $\#(E_y\cap Y)\leq N$, so $Y$ has bounded geometry, proving $(i)\implies (iii)$.

Suppose $Y\subseteq X$ is a coarsely dense discrete space of bounded geometry and let $\mu$ be the counting measure of $Y$.
Any controlled set $E$ on $X$ restricts to a controlled set on $Y$.
Since $Y$ is a discrete space of bounded geometry any $E$-bounded set can contain a bounded number of points of $Y$.
This shows that $\mu$ is uniformly bounded.
Now let $E\subseteq X\times X$ be a symmetric controlled set such that for all $x\in X$ there is $y\in Y$ with $(x,y)\in E$.
Then each $E_x$ has measure at least 1, so $E$ is gordo, showing $(iii)\implies (ii)$.

Finally, let $\mu$ be a uniformly bounded measure on $X$ and let $F$ be a gordo set for $\mu$.
We will show that $F^{\circ4}$ is a covering set.
There is $\epsilon>0$ such that $\mu(F_x)\geq\epsilon$ for all $x$.
Let $E$ be a symmetric controlled set.
Since $\mu$ is uniformly bounded, there is a constant $C$ such that $\mu(V)\leq C$ for all $F\circ E\circ F$-bounded measurable $V$.
Now let $U\subseteq X$ be $E$-bounded.
Then $F_U$ is $F\circ E\circ F$-bounded, so all its measurable subsets have measure at most $C$.
Let $Y\subseteq U$ be maximal such that the $F_y$ are pairwise disjoint when the $y$ range over $Y$.
Since the $F_y$ are all contained in $F_U$ and have measure at least $\epsilon$ we have $\#Y\leq N$, with $N=\frac C\epsilon$.
Since $Y$ is maximal we know that for all $x\in U$ there is a $y\in Y$ with $F_x\cap F_y\neq\emptyset$, so $(x,y)\in F\circ F$.
This implies that the sets $(F\circ F)_y,y\in Y$ cover $U$.
These sets are $F^{\circ4}$-bounded.
By Lemma \ref{lem-different_def_covering} we see that $F^{\circ4}$ is covering, showing $(ii)\implies(i)$.
\end{proof}

\section{Blocking collections}

From now on, let $X$ be a coarse space of bounded geometry that can be covered by a countable number of bounded sets.
Let $\mu$ be a uniformly bounded measure on $X$ for which a gordo set exists.
This measure exists by Proposition \ref{prop-bg-coarse}; in applications, there is often already a measure given that satisfies the condition above.
For example, if $X$ is a metric space and $\mu$ is a Borel measure, then $\mu$ satisfies the conditions if there is a radius $R>0$ such that $R$-balls are bounded from below in measure, and for all radii $R>0$ the $R$-balls are bounded from above in measure.

The definition of geometric property (T) uses the measure $\mu$, but we will show later that it does not depend on it (see Theorem \ref{thm-coarse-invariance}).

We will start by defining the notion of a blocking collection of subsets of $X$ and proving some elementary results about them.

\begin{definition}
Let $E$ be a controlled set and let $\epsilon>0$.
A collection $(A_i)$ of measurable subsets of $X$ is called \emph{$(E,\epsilon)$-blocking} if $\bigcup_iA_i\times A_i\subseteq E$, the $A_i$ are pairwise disjoint and $\mu(A_i)\geq\epsilon$ for all $i$.
A collection $(A_i)$ of measurable subsets of $X$ is called \emph{blocking} if it is $(E,\epsilon)$-blocking for some $E,\epsilon$.
\end{definition}

The conditions on the measure $\mu$ ensure that there are sufficiently large blocking collections:
\begin{lemma}\label{lem-complete_blocking_collection}
There is a blocking collection $(A_i)$ with union $X$.
\end{lemma}
\begin{proof}
Let $E\subseteq X\times X$ be a gordo set for $\mu$.
By Zorn's Lemma, there is a maximal subset $I\subseteq X$ such that all the $E_i$ are pairwise disjoint when $i$ ranges over $I$.
Then for all $x\in X$ there exists an $i\in I$ such that $E_x\cap E_i\neq\emptyset$, so $(x,i)\in E\circ E$.
So the $(E\circ E)_i$ cover $X$, while the $E_i$ are pairwise disjoint.
Hence, we can find measurable $A_i$ with $E_i\subseteq A_i\subseteq (E\circ E)_i$ such that $X=\bigsqcup_iA_i$.
To construct the $A_i$ explicitly, choose a well-ordering on $I$ and define
\[A_i=E_i\cup\left((E\circ E)_i\setminus\left(\bigcup_{i'\neq i}E_{i' }\cup \bigcup_{i'<i}(E\circ E)_{i'}\right)\right).\]
Now the measure of the $A_i$ is bounded from below because $E$ is gordo, and the $A_i$ are $E^{\circ4}$-bounded, showing that the $A_i$ form a blocking collection.
\end{proof}

\begin{definition}\label{def-controlled-set-in-blocks}
Let $E\subseteq X\times X$.
We say that $E$ is a \emph{controlled set in blocks} if $E$ is of the form $\bigsqcup_iA_i\times A_i$, where the $(A_i)$ form a blocking collection.
\end{definition}
Note that $\bigsqcup_iA_i\times A_i$ is always a controlled set if the $(A_i)$ form a blocking collection.

\begin{lemma}\label{lem-union-of-sets-in-blocks}
Let $E\subseteq X\times X$ be a controlled set.
Then $E$ is contained in a finite union of controlled sets in blocks.
\end{lemma}
\begin{proof}
Assume without loss of generality that $E$ is symmetric.
By Lemma \ref{lem-complete_blocking_collection}, there is a blocking collection $(A_i)$ with union $X$.
Let $F=\bigsqcup_iA_i\times A_i$ which is a controlled set.
Define a graph structure on the index set $I$ by connecting $i$ and $j$ if $E_{A_i}\cap E_{A_j}\neq\emptyset$.

If $i$ and $j$ are connected by an edge, then $(E\circ E)_{A_i}\cap A_j\neq \emptyset$, and since $A_j$ is $F$-bounded it follows that $A_j\subseteq (F\circ E\circ E)_{A_i}$.
Since $\mu$ is uniformly bounded and the $A_i$ are $F$-bounded, we know that the measure of $(F\circ E\circ E)_{A_i}$ is bounded from above, say by $M$, where $M$ does not depend on $i$.
Moreover, the measure of $A_j$ is bounded from below uniformly in $j$, say by $\epsilon>0$.
Since the $A_j$ are also pairwise disjoint, the degree of any vertex $i$ in the graph can be at most $\frac M\epsilon$.
Then we can colour the graph in $N$ colours, for some integer $N$, such that no two adjacent vertices have the same colour.
Now for $1\leq k\leq N$ define
\[\mathcal C_k=\{E_{A_i}\mid i\text{ is coloured with the }k\text{-th colour}\}.\]
Then because of the way we coloured the graph, the sets in $\mathcal C_k$ are disjoint.
Moreover they are $E\circ F\circ E$-bounded and bounded from below in measure, so $\mathcal C_k$ is a blocking collection for each $k$.
For $1\leq k\leq N$ we now define
\[E_k=\bigsqcup\{E_{A_i}\times E_{A_i}\mid i\text{ is coloured with the }k\text{-th colour}\}.\]
These are controlled sets in blocks.
For all $(x,y)\in E$, there is an $i$ with $y\in A_i$.
Let $k$ be the colour of $i$.
Then $(x,y)\in E_{A_i}\times E_{A_i}\subseteq E_k$.
So $E$ is contained in $\bigsqcup_kE_k$.
\end{proof}

\section{The algebra of a coarse space}\label{sec-algebra-coarse-space}

Analogous to the discrete setting, we now define the algebra of operators with controlled support.
We will use the Hilbert space $L^2X=L^2(X,\mu)$.
\begin{definition}\label{def-controlled-support}
Let $E\subseteq X\times X$ be a measurable controlled set and $T\in B(L^2X)$.
We say that $T$ is \emph{supported on $E$} if for all measurable $U\subseteq X$, all $\xi\in L^2X$ supported on $U$ and all $\eta\in L^2X$ we have
\[\langle\xi,T\eta\rangle=\langle\xi,T\eta|_{E_U}\rangle.\]
We say that $T$ has \emph{controlled support} if it is supported on some measurable controlled set.
\end{definition}

We can give some equivalent definitions for this:
\begin{lemma}\label{lem-eq-defs-controlled-support}
Let $E$ be a measurable symmetric controlled set and $T\in B(L^2X)$.
The following are equivalent:
\begin{enumerate}[(i)]
    \item The operator $T$ is supported on $E$.
    \item For all measurable $U\subseteq X$, all $\xi\in L^2X$ and all $\eta\in L^2X$ supported on $U$ we have
    \[\langle\xi,T\eta\rangle=\langle\xi|_{E_U},T\eta\rangle.\]
    \item The operator $T^*$ is supported on $E$.
    \item For all measurable $U\subseteq X$ and all $\eta\in L^2X$ supported on $U$ the function $T\eta$ is supported on $E_U$.
\end{enumerate}
\end{lemma}
\begin{proof}
To prove $(i)\implies(ii)$, we have
\[\langle\xi,T\eta\rangle=\langle\xi|_{E_U},T\eta\rangle+\langle\xi|_{X\setminus E_U},T\eta\rangle=\langle\xi|_{E_U},T\eta\rangle+\langle\xi|_{X\setminus E_U},T\eta|_{E_{X\setminus E_U}}\rangle=\langle\xi|_{E_U},T\eta\rangle.\]
The other direction is similar, and the equivalences $(ii)\iff (iii)$ and $(ii)\iff (iv)$ are clear.
\end{proof}

Clearly linear combinations of operators with controlled support again have controlled support, as does a product of operators with controlled support.
So the operators with controlled support form a unital $*$-algebra.
\begin{definition}\label{def-algebra-controlled-support}
Denote by $\mathbb C_{\cs}[X]$ the subalgebra $\{T\in B(L^2X)\mid T\text{ has controlled support}\}$ of $B(L^2X)$.
We will denote it $\mathbb C_{\cs}[X,\mu]$ when the measure is not clear from the context.
\end{definition}
\begin{remark}
In the case that $X$ is a discrete space of bounded geometry, the algebra $\mathbb C_{\cs}[X]$ is the same as the algebra $\mathbb C_u[X]$ as defined in \cite[Definition 3.1]{WilYu14}.
\end{remark}

\begin{definition}\label{def-operators-in-blocks}
Let $T\in\mathbb C_{\cs}[X]$, let $E$ be a controlled set and let $\epsilon>0$.
We say that $T$ is an \emph{operator in $(E,\epsilon)$-blocks} if it is supported on $\bigsqcup_iA_i\times A_i$ where $(A_i)$ is some $(E,\epsilon)$-blocking collection.
We say that $T$ is an \emph{operator in blocks} if it is supported on some controlled set in blocks.
\end{definition}

For a measurable subset $U\subseteq X$, denote by $\pi_U\colon L^2X\to L^2U$ the restriction map and by $i_U\colon L^2U\to L^2X$ extension by 0.
An operator $T$ that is supported on $\bigsqcup_iA_i\times A_i$ can be written as $T=\sum_iT_i$, where $T_i=\pi_{A_i}Ti_{A_i}\colon L^2A_i\to L^2A_i$, using Lemma \ref{lem-eq-defs-controlled-support}.
This decomposition makes it very convenient to do computations with operators in blocks.
Moreover, every operator in $\mathbb C_{\cs}[X]$ can be written as a finite sum of operators in blocks:
\begin{lemma}\label{lem-sum-of-blocks}
Let $T\in\mathbb C_{\cs}[X]$.
We can write $T$ as a finite sum $T_1+\ldots+T_N$ where each $T_i\in \mathbb C_{\cs}[X]$ is an operator in blocks.
\end{lemma}
\begin{proof}
Let $T$ be supported on the measurable and symmetric controlled set $E$.
Let $(A_i)$ be a blocking collection as in the proof of Lemma \ref{lem-union-of-sets-in-blocks} and consider the same colouring on $I$ as in that proof.
Let
\[T_k=\sum_{i\text{ with colour }k}\pi_{A_i}T.\]
Then $T_k$ is controlled on $\sqcup_{i\text{ with colour }k}E_{A_i}\times E_{A_i}$, so it is an operator in blocks.
Moreover we have $T=T_1+\ldots+T_N$.
\end{proof}

We define some spaces of functions on $X$:
\begin{definition}
\begin{enumerate}[(i)]
    \item Let $L^2_bX$ denote the $L^2$-functions on $X$ with bounded support.
    \item Let $LL^2X$ denote the measurable functions on $X$ that are locally $L^2$: for each $\phi\in LL^2X$ and each bounded measurable $U$, we have $\int_U|\phi|^2<\infty$.
\end{enumerate}
\end{definition}

\begin{lemma}
The dual of $L^2_bX$ can be identified with $LL^2X$.
\end{lemma}
\begin{proof}
For any $\xi\in L^2_bX$ and $\eta\in LL^2X$, the inner product $\langle\eta,\xi\rangle$ is finite, so $\eta$ defines a function in $(L^2_bX)^*$.

Conversely, let $\phi\in(L^2_bX)^*$.
For any bounded measurable $U$, the function $\phi$ can be pulled back along the inclusion $L^2U\hookrightarrow L^2_bX$ to an element of $(L^2U)^*$, which may be viewed as a function $\eta_U\in L^2U$.
If $U\subseteq V$ are bounded and measurable, then we find that $\eta_U=\eta_V|_{U}$.
Therefore, the $\eta_U$ may be glued together to a function $\eta\in LL^2X$.
\end{proof}

If $T\in B(L^2X)$ has controlled support, it sends $L^2_bX$ to $L^2_bX$, and so does the dual $T^*$.
Taking the dual of $T^*\colon L^2_b(X)\to L^2_b(X)$ gives a map $LL^2X\to LL^2X$ which we denote by $T$ again.
Now define the linear map $\Phi\colon \mathbb C_{\cs}[X]\to LL^2X$ by $\Phi(T)=T(\mathbb 1_X)\in LL^2X$.
Explicitly, if $T$ is supported on a controlled set $E$ then we have
\[\Phi(T)(x)=(T\mathbb 1_{E_x})(x).\]
For all $T,S\in \mathbb C_{\cs}[X]$ we have $\Phi(ST)=S\Phi(T)$.

\begin{lemma}\label{lem-sum-of-blocks-with-Phi-bounded}
Let $T\in\mathbb C_{\cs}[X]$ be such that $\Phi(T)$ is bounded.
Then we can write $T=T_1+\ldots+T_N$ where the $T_k$ are operators in blocks and $\Phi(T_k)$ is bounded.
\end{lemma}
\begin{proof}
Construct $T_k$ as in the proof of Lemma \ref{lem-sum-of-blocks}.
Then the $\Phi(T_k)$ are bounded.
\end{proof}

\section{Geometric property (T)}

\begin{definition}
A representation of $\mathbb C_{\cs}[X]$ consists of a Hilbert space $\H$ and a unital $*$-homomorphism $\rho\colon \mathbb C_{\cs}[X]\to B(\H)$.
If no confusion arises we will omit the map $\rho$ from the notation.
\end{definition}

The standard representation of $\mathbb C_{\cs}[X]$ is the inclusion $\mathbb C_{\cs}[X]\hookrightarrow B(L^2X)$.
Note that we do not ask the map $\rho$ to be continuous with respect to the norm on $\mathbb C_{\cs}[X]$, inherited from $B(L^2X)$.

\begin{definition}
Let $T,S\in\mathbb C_{\cs}[X]$.
\begin{enumerate}[(i)]
    \item We write $T\geq S$ if this inequality holds in $B(L^2X)$.
    \item We write $T\geq_{\max} S$ if for every representation $(\rho,\H)$ of $\mathbb C_{\cs}[X]$, we have the inequality $\rho(T)\geq\rho(S)$.
    \item Let $\norm T_{L^2}$ be the norm of $T$ in $B(L^2X)$.
    \item Define $\norm T_{\max}=\sup\norm{\rho(T)}$, where the supremum is over all representations $(\rho,\H)$.
    We will see in Corollary \ref{cor-max-norm-finite} that this is always finite, and hence it also defines a norm on $\mathbb C_{\cs}[X]$.
\end{enumerate}
\end{definition}

\begin{lemma}\label{lem-maximal-things-for-operators-in-blocks}
Let $T\in\mathbb C_{\cs}[X]$ be an operator in blocks.
\begin{enumerate}[(i)]
    \item If $T\geq 0$, then $T\geq_{\max} 0$.
    \item We have $\norm T_{L^2}=\norm T_{\max}$.
\end{enumerate}
\end{lemma}
\begin{proof}
\begin{enumerate}[(i)]
    \item Suppose $T$ is supported on $\bigsqcup_iA_i\times A_i$ where $(A_i)$ is a blocking collection.
    Write $T=\sum_iT_i$ with $T_i\colon L^2A_i\to L^2A_i$.
    Then $T^\frac12=\sum_iT_i^\frac12$ is again an operator in blocks.
    Now for any representation $(\rho,\H)$ we have $\rho(T)=\rho(T^\frac12)^2\geq 0$, so we have $T\geq_{\max}0$.
    \item Clearly we have $\norm T_{L^2}\leq\norm T_{\max}$.
    Note that $T^*T$ is again an operator in blocks.
    We have $0\leq T^*T\leq \norm T_{L^2}^2$, and by part (i), it follows that $0\leq_{\max} T^*T\leq_{\max} \norm T_{L^2}^2$.
    So for any representation $(\rho,\H)$, we have $0\leq \rho (T^*T)\leq \norm T_{L^2}^2$.
    So $\norm{\rho(T)}=\norm{\rho(T)^*\rho(T)}^\frac12=\norm{\rho(T^*T)}^\frac12\leq \norm T_{L^2}$, and hence $\norm T_{L^2}=\norm T_{\max}$.
\end{enumerate}
\end{proof}
If $T$ is an operator in blocks we will simply write $\norm T$ for both the norms $\norm T_{L^2}$ and $\norm T_{\max}$.
\begin{corollary}\label{cor-max-norm-finite}
For any $T\in\mathbb C_{\cs}[X]$, the maximal norm $\norm T_{\max}$ is finite.
\end{corollary}
\begin{proof}
This follows directly from Lemma \ref{lem-maximal-things-for-operators-in-blocks} and Lemma \ref{lem-sum-of-blocks}.
\end{proof}
Therefore $\norm{\cdot}_{\max}$ defines a norm and we can define the maximal $C^*$-algebra:
\begin{definition}\label{def-maximal-C^*-algebra}
The maximal $C^*$-algebra $C^*_{\max}(X)$ is defined as the completion of $\mathbb C_{\cs}[X]$ with respect to the maximal norm $\norm\cdot_{\max}$.
\end{definition}

For a discrete space $X$ of bounded geometry the $C^*$-algebra $C^*_{\max}(X)$ is the same as the one defined in \cite{WilYu14} (there it is called $C^*_{u,\max}(X)$).
For $f\in L^\infty X$ denote by $M_f\in B(L^2X)$ the multiplication operator by $f$.
Similar to \cite{WilYu14}, we define the subspace of the constant vectors.
\begin{definition}\label{def-constant-vector}
Let $\H$ be a non-degenerate representation of $\mathbb C_{\cs}[X]$, and let $v\in\H$.
We call $v$ a \emph{constant vector} if $Tv=M_{\Phi(T)}v$ for all operators $T\in \mathbb C_{\cs}[X]$ for which $\Phi(T)$ is bounded.
Here $\Phi$ is as defined at the end of Section \ref{sec-algebra-coarse-space}.
The constant vectors form a closed subspace $\H_c$ of $\H$.
\end{definition}

Since we can replace $T$ by $T-M_{\Phi(T)}$, it is enough to ask that $Tv=0$ for all $T$ with $\Phi(T)=0$.

Now we are ready to define geometric property (T).

\begin{definition}\label{def-geometric-T}
Let $X$ be a coarse space with bounded geometry and let $\mu$ be a uniformly bounded measure for which a gordo set exists.
We say that $(X,\mu)$ satisfies \emph{geometric property (T)} if there are a controlled set $E$ and constants $\epsilon,\gamma>0$ such that for every non-degenerate representation $(\rho,\H)$ and every unit vector $v\in\H_c^{\perp}$ there is an operator $T$ in $(E,\epsilon)$-blocks with $\Phi(T)$ bounded and $\norm{(T-M_{\Phi(T)})v}\geq \gamma\norm T$.
\end{definition}

In the definition, it is necessary to ask that $\Phi(T)$ is a bounded function, because otherwise $M_{\Phi(T)}$ would not be defined.
By considering the operator $T-M_{\Phi(T)}$, we see that we can equivalently ask there to be an operator $T$ in $(E,\epsilon)$-blocks with $\Phi(T)=0$ and $\norm{Tv}\geq \gamma\norm T$.

If $X$ is a discrete space of bounded geometry with a generating controlled set, and $\mu$ is the counting measure, we retrieve the definition of \cite{WilYu14} (cf. Proposition 3.1 in that paper):

\begin{lemma}\label{lem-(T)-the-same-as-in-WilYu}
Let $X$ be a discrete space of bounded geometry for which there is a generating controlled set, and let $\mu$ be the counting measure.
The following are equivalent: 
\begin{enumerate}[(i)]
    \item The space $(X,\mu)$ satisfies geometric property (T).
    \item There exists a controlled generating set $E$ and a constant $\gamma>0$ such that for any representation $\H$ and unit vector $v\in\H_c^{\perp}$, there is an operator $T\in\mathbb C_{\cs}[X]$ with support in $E$ such that
    \[\norm{(T-M_{\Phi(T)})v}> \gamma\sup_{x,z}|T_{xz}|,\]
    where $T_{xz}$ denotes the value at $(x,z)$ of the matrix corresponding to the operator $T$.
\end{enumerate}
\end{lemma}
\begin{proof}
We certainly have $\norm{T}_{L^2}\geq \sup_{x,z}|T_{xz}|$.
With this, the direction $(i)\implies(ii)$ is trivial, since we may enlarge the controlled set $E$ to make it generating.

Now suppose $(ii)$ is true, take $E,\gamma$ satisfying the condition, and let $\epsilon=1$.
Without loss of generality, we may assume $E$ to be symmetric.
There is an integer $N$ such that $\#E_x\leq N$ for each $x\in X$.
Let $\H$ be a representation and $v\in \H_c^{\perp}$ a unit vector.
We know there is $T\in\mathbb C_{\cs}[X]$ supported on $E$ such that $\norm{(T-\Phi(T))v}>\gamma \sup_{x,z}|T_{xz}|$.
For each pair $(x,z)\in E$ there are fewer than $4N$ other pairs $(x',z')\in E$ with $\{x,z\}\cap\{x',z'\}\neq\emptyset$.
By a graph colouring argument we may write $T=T_1+\ldots+T_{4N}$ as a sum of operators controlled on $E$ with for each $k$: if $(x,z)\in E$ and $(x',z')\in E$ are different pairs where $T_k$ is supported, then $\{x,z\}\cap\{x',z'\}=\emptyset$.
Moreover we get $\sup_{x,z}|T_{xz}|=\sup_k\norm{T_k}_{L^2}$.
We then have
\[\norm{\sum_k(T_k-M_{\Phi(T_k)})v}=\norm{(T-M_{\Phi(T)})v}\geq \gamma\sup_{x,z}|T_{xz}|=\gamma\sup_k\norm{T_k}_{L^2},\]
hence there is a $k$ with $\norm{(T_k-M_{\Phi(T_k)})v}\geq\frac\gamma{4N}\norm{T_k}_{L^2}$, proving $(i)$.
\end{proof}

\section{Laplacians}

In this section we will give a different characterisation of geometric property (T) using spectral properties of Laplacians.

\begin{definition}\label{def-laplacian-operator}
Let $E$ be a measurable symmetric controlled set.
Define the \emph{Laplacian} $\Delta^E\in B(L^2X)$ as 
\[\Delta^E(\xi)(x)=\int_{E_x}(\xi(x)-\xi(y))d\mu(y).\]
\end{definition}

The Laplacian corresponding to $E$ is obviously supported on $E$, in particular we have $\Delta^E\in\mathbb C_{\cs}[X]$.
The Laplacians are positive operators, as is shown by the following lemma.
\begin{lemma}\label{lem-laplacians-positive-and-bounded}
For any measurable symmetric controlled set $E$ we have $0\leq_{\max}\Delta^E\leq_{\max} 2M$,
where $M=\sup_x\mu(E_x)$.
\end{lemma}
\begin{proof}
A straightforward computation shows that for all measurable symmetric controlled sets $F$, we have
\begin{align*}
    \langle\Delta^F\xi,\xi\rangle&=\int_X\overline{\Delta^F\xi(x)}\xi(x)d\mu(x)\\
    &=\int_X\int_{F_x}\overline{(\xi(x)-\xi(y))}\xi(x)d\mu(y)d\mu(x)\\
    &=\frac12\int_X\int_{F_x}\left(\overline{(\xi(x)-\xi(y))}\xi(x)+\overline{(\xi(y)-\xi(x))}\xi(y)\right)d\mu(y)d\mu(x)\\
    &=\frac12\int_F|\xi(x)-\xi(y)|^2d(\mu\times\mu)(x,y)\\
    &\geq0,
\end{align*}
so $\Delta^F\geq 0$.
Since $|\xi(x)-\xi(y)|^2\leq 2|\xi(x)|^2+2|\xi(y)|^2$ we also have
\[\langle\Delta^F\xi,\xi\rangle\leq\int_F(|\xi(x)|^2+|\xi(y)|^2)d(\mu\times\mu)(x,y)=2\int_F|\xi(x)|^2d(\mu\times\mu)(x,y)\]
by symmetry.
This is equal to $2\int_F|\xi(x)|^2\mu(F_x)d\mu(x)$, which is at most $2\norm\xi^2\sup_x\mu(F_x)$, so we have $\norm{\Delta^F}_{L^2}\leq 2\sup_x\mu(F_x)$.
Moreover, if $F_1,F_2$ are disjoint measurable symmetric controlled sets, we see directly that $\Delta^{F_1\cup F_2}=\Delta^{F_1}+\Delta^{F_2}$.

Let $E$ be any measurable symmetric controlled set.
By Lemma \ref{lem-union-of-sets-in-blocks}, there are symmetric controlled sets in blocks $E_1,\ldots,E_N$ whose union contains $E$.
Let $E'_1,\ldots,E'_N$ be pairwise disjoint measurable sets such that $E'_k\subseteq E_k$ for each $k$ and $E=\bigsqcup_kE'_k$.
For each $k$, the operator $\Delta^{E'_k}$ is an operator in blocks, so by Lemma \ref{lem-maximal-things-for-operators-in-blocks} we have $0\leq_{\max}\Delta^{E'_k}\leq_{\max}2\sup_x\mu((E'_k)_x)$.
Taking the sum from $k=1$ to $N$, we get $0\leq_{\max}\Delta^E\leq_{\max} 2\sup_x\mu(E_x)$.
\end{proof}

\begin{lemma}\label{lem-bounded-above-by-laplacian}
Let $E$ be a controlled set in blocks and let $T\in\mathbb C_{\cs}[X]$ be a positive operator supported on $E$ such that $\Phi(T)=0$.
Then
\[T\leq_{\max}\frac{\norm T}{\inf_x(\mu(E_x))}\Delta^E,\]
where the infimum is taken over all $x$ with $E_x\neq\emptyset$.
\end{lemma}
\begin{proof}
Write $E=\bigsqcup_iA_i\times A_i$ and $T=\sum_iT_i$ with $T_i\colon L^2A_i\to L^2A_i$.
It is enough to prove the inequality $T_i\leq\frac{\norm{T_i}}{\mu(A_i)}\Delta^{E_i}$ for each $i$, where $E_i=A_i\times A_i$.
The space $L^2A_i$ can be written as the direct sum of the constant functions and the functions with zero integral.
For any $\xi\in L^2A_i$ write $\xi=\xi_c+\xi_d$ where $\xi_c$ is constant and $\int\xi_d=0$.
Then $T_i\xi=T_i\xi_d$ and $\Delta^{E_i}\xi=\mu(A_i)\xi_d$.
So we have $\langle T_i\xi,\xi\rangle\leq\norm{T_i}\cdot\norm{\xi_d}^2=\langle \frac{\norm{T_i}}{\mu(A_i)}\Delta^{E_i}\xi,\xi\rangle$ showing that $T_i\leq\frac{\norm{T_i}}{\mu(A_i)}\Delta^{E_i}$.
\end{proof}

Since $\Phi(\Delta^E)=0$, we have $\Delta^Ev=0$ for all $v\in\H_c$, where $\H$ is a representation of $\mathbb C_{\cs}[X]$.
Moreover, we have the following lemma, which shows the converse.
\begin{lemma}\label{lem-intersection-kernels-laplacians}
Let $(\rho,\H)$ be a representation of $\mathbb C_{\cs}[X]$.
Then
\[\H_c=\cap_E\ker(\rho(\Delta^E))\]
where the intersection ranges over all the measurable symmetric controlled sets $E$.
\end{lemma}
\begin{proof}
Let $v\in\cap_E\ker(\rho(\Delta^E))$ and let $T\in\mathbb C_{\cs}[X]$ such that $\Phi(T)$ is bounded, we will show that $Tv=M_{\Phi(T)}v$, showing that $v\in \H_c$.
By Lemma \ref{lem-sum-of-blocks-with-Phi-bounded}, we can assume that $T$ is an operator in blocks.
We can also assume that $\Phi(T)=0$ by considering the operator $T-\Phi(T)$.
By Lemma \ref{lem-bounded-above-by-laplacian}, there is a measurable symmetric controlled set $E$ and a constant $c$ such that $T^*T\leq_{\max} c\Delta^E$, since $T^*T$ is a positive operator in blocks.
Then $\Delta^Ev=0$ and we get $\norm{Tv}^2=\langle T^*Tv,v\rangle\leq c\langle \Delta^Ev,v\rangle=0$, so $Tv=0$.
\end{proof}

\begin{definition}\label{def-spectral-gap}
Let $T\in C^*_{\max}(X)$ be a positive operator.
We say that $T$ \emph{has spectral gap} if there is a constant $\gamma>0$ such that $\sigma_{\max}(T)\subseteq\{0\}\cup[\gamma,\infty)$, where $\sigma_{\max}$ denotes the spectrum of $T$ in the maximal $C^*$-algebra.
\end{definition}

Now we can prove the generalisation of \cite[Proposition 5.2]{WilYu14}.

\begin{proposition}\label{prop-property-(T)-iff-laplacian-spectral-gap}
Let $X$ be a coarse space of bounded geometry and let $\mu$ be a uniformly bounded measure for which a gordo set exists.
Then $(X,\mu)$ has geometric property (T) if and only if there exists a measurable symmetric controlled set $E$ such that for each representation $(\rho,\H)$ we have $\H_c=\ker(\rho(\Delta^E))$ and $\Delta^E$ has spectral gap.
\end{proposition}
\begin{proof}
Suppose $X$ has property (T).
Let $E$ be a controlled set and $\epsilon,\gamma>0$ constants such that for every representation $\H$ and every unit vector $v\in\H_c^\perp$, there is an operator $T$ in $(E,\epsilon)$-blocks such that $\Phi(T)=0$ and $\norm {Tv}\geq \gamma\norm T$.
We will show that for each representation $(\rho,\H)$, we have $\H_c=\ker(\rho(\Delta^E))$ and $\Delta^E$ has spectral gap.

Let $v\in\H_c^\perp$ be a unit vector.
Let $T$ be an operator in $(E,\epsilon)$-blocks with $\Phi(T)=0$ and $\norm{Tv}\geq \gamma\norm T$.
Let $T$ be supported on the blocks $(A_i)$ and let $E'=\bigsqcup_iA_i\times A_i\subseteq E$.
By Lemma \ref{lem-bounded-above-by-laplacian}, we have $T^*T\leq \frac{\norm T^2}{\inf_i\mu(A_i)}\Delta^{E'}\leq\frac{\norm T^2}\epsilon\Delta^E$.
Now
\[\gamma^2\norm T^2\leq \norm{Tv}^2=\langle T^*Tv,v\rangle\leq \frac{\norm T^2}\epsilon\langle\Delta^Ev,v\rangle\leq \frac{\norm T^2}\epsilon\norm{\Delta^Ev},\]
so $\norm{\Delta^Ev}\geq \gamma^2\epsilon$.
This shows that $\ker(\rho(\Delta^E))\cap\H_c^\perp=\{0\}$, so $\ker(\rho(\Delta^E))=\H_c$.
Moreover, it shows that $\sigma(\rho(\Delta^E))\subseteq \{0\}\cup[\gamma^2\epsilon,\infty)$.
Since this holds for every representation $(\rho,\H)$, we get $\sigma_{\max}(\Delta^E)\subseteq\{0\}\cup[\gamma^2\epsilon,\infty)$, hence $\Delta^E$ has spectral gap.

Now suppose that $E$ is a measurable symmetric controlled set such that for each representation $(\rho,\H)$ we have $\H_c=\ker(\rho(\Delta^E))$ and $\Delta^E$ has spectral gap, say $\sigma_{\max}(\Delta^E)\subseteq\{0\}\cup[\gamma,\infty)$.
By Lemma \ref{lem-union-of-sets-in-blocks} there are controlled sets in blocks $E_1,\ldots,E_N$ whose union contains $E$.
Let $\epsilon=\min_k\inf_x\mu((E_k)_x)>0$, where $k$ ranges from 1 to $N$ and $x$ ranges over all elements in $X$ with $(E_k)_x\neq\emptyset$.
Now let $(\rho,\H)$ be a representation and let $v\in\H_c^\perp$ be a unit vector.
Then $v\in\ker(\rho(\Delta^E))^\perp$ so $\langle\Delta^Ev,v\rangle\geq\gamma$.
We also have $\Delta^E\leq\Delta^{E_1}+\ldots+\Delta^{E_N}$ so there is $1\leq k\leq N$ such that $\langle\Delta^{E_k}v,v\rangle\geq\frac\gamma N$.
Now $\Delta^{E_k}$ is an $(E,\epsilon)$-operator in blocks and $\norm{\Delta^{E_k}v}\geq\frac\gamma N$, so we see that $X$ has geometric property (T).
\end{proof}

In \cite[Proposition 5.2]{WilYu14} it is actually shown that $\H_c=\ker(\rho(\Delta^E))$ holds for all generating symmetric controlled sets $E$ and that if $X$ has geometric property (T), then $\Delta^E$ has spectral gap for all generating symmetric controlled sets $E$.
In our case, we have a similar result, for which we first need another lemma and its corollary.

\begin{lemma}\label{lem-bound-laplacian-composition}
Let $E$ be a gordo set and $F$ a symmetric measurable controlled set containing $E$ and satisfying $F\circ E=E\circ F$.
Let $f\in L^\infty(X)$ be the function $f(x)=\mu((E\circ F)_x)$.
Then there is a constant $\delta$ such that
\[\delta\Delta^{F\circ F}\leq_{\max}(2M_f-\Delta^{E\circ F})M_{\frac1f}\Delta^{E\circ F}.\]
\end{lemma}
\begin{proof}
For every $\xi\in L^2X$ and $x\in X$, we have
\[(2M_f-\Delta^{E\circ F}
)\xi(x)=f(x)\xi(x)+\int_{(E\circ F)_x}\xi(y)d\mu(y).\]
Then we have
\begin{align*}
&(2M_f-\Delta^{E\circ F})M_{\frac1f}\Delta^{E\circ F}\xi(x)\\
=&\Delta^{E\circ F}\xi(x)+\int_{(E\circ F)_x}\frac1{f(y)}\Delta^{E\circ F}\xi(y)d\mu(y)\\
=&f(x)\xi(x)-\int_{(E\circ F)_x}\xi(y)d\mu(y)+\int_{(E\circ F)_x}\xi(y)d\mu(y)-\int_{(E\circ F)_x}\frac1{f(y)}\int_{(E\circ F)_y}\xi(z)d\mu(z)d\mu(y)\\
=&f(x)\xi(x)-\int_{(E\circ F)^{\circ2}_x}\alpha(x,z)\xi(z)d\mu(z),
\end{align*}
where
\[\alpha(x,z)=\int_{(E\circ F)_x\cap(E\circ F)_z}\frac1{f(y)}d\mu(y).\]
Let $z\in (F\circ F)_x$.
There is $y\in X$ with $(x,y)\in F$ and $(y,z)\in F$.
It follows that $E_y\subseteq (E\circ F)_x\cap (E\circ F)_z$, so $\alpha(x,z)\geq \int_{E_y}\frac1{f(y')}d\mu(y')$.
Since $E$ is gordo and $f$ is bounded from above, there is a constant $\delta$ such that $\alpha(x,z)\geq \delta$ for all $z\in (F\circ F)_x$.
Now define
\[\beta(x,z)=\begin{cases}\alpha(x,z)-\delta&\text{ if }(x,z)\in F\circ F,\\\alpha(x,z)&\text{ otherwise.}\end{cases}\]
Then $\beta$ is symmetric, $\beta(x,z)\geq 0$ for all $(x,z)\in X\times X$ and
\[(2M_f-\Delta^{E\circ F})M_{\frac1f}\Delta^{E\circ F}=\delta\Delta^{F\circ F}+T,\]
where $T\in\mathbb C_{cs}[X]$ is the operator defined by
\[T\xi(x)=\int_{(E\circ F)^{\circ2}_x}\beta(x,y)(\xi(x)-\xi(y))d\mu(y).\]
It remains to show that $T\geq_{\max} 0$.
This is proved similarly to Lemma \ref{lem-laplacians-positive-and-bounded}.
Let $E_1,\ldots,E_N$ be controlled sets in blocks whose union contains $(E\circ F)^{\circ2}$.
Let $E_1',\ldots,E_N'$ be measurable symmetric subsets of $E_1,\ldots,E_N$ respectively, that are disjoint from each other and whose union is $(E\circ F)^{\circ2}$.
For $1\leq k\leq N$ define $T_k\in\mathbb C_{\cs}[X]$ by
\[T_k\xi(x)=\int_{(E'_k)_x}\beta(x,y)(\xi(x)-\xi(y))d\mu(y).\]
An easy computation shows that
\[\langle T_k\xi,\xi\rangle=\frac12\int_{E'_k}\beta(x,y)(\xi(x)-\xi(y))^2d(\mu\times\mu)(x,y),\]
so $T_k\geq 0$.
Since $T_k$ is an operator in blocks, we get $T_k\geq_{\max}0$ by Lemma \ref{lem-maximal-things-for-operators-in-blocks}.
So we have $T=\sum_{k=1}^NT_k\geq_{\max}0$.
\end{proof}

\begin{corollary}\label{cor-laplacians-of-close-sets-behave-the-same}
Let $E$ be a gordo set and $F$ a symmetric measurable controlled set satisfying $F\circ E=E\circ F$ and $E\subseteq F$.
Then for every representation $(\rho,\H)$, we have $\ker(\rho(\Delta^{F\circ F}))=\ker(\rho(\Delta^{E\circ F}))$.
Moreover, $\Delta^{F\circ F}$ has spectral gap if and only if $\Delta^{E\circ F}$ has spectral gap.
\end{corollary}
\begin{proof}
Since $E\subseteq F$ we have $E\circ F\subseteq F\circ F$, and therefore $\Delta^{E\circ F}\leq_{\max}\Delta^{F\circ F}$.
It follows that $\ker(\rho(\Delta^{F\circ F}))\subseteq\ker(\rho(\Delta^{E\circ F}))$ and that if $\Delta^{E\circ F}$ has spectral gap, then $\Delta^{F\circ F}$ also has spectral gap.
The other direction follows from Lemma \ref{lem-bound-laplacian-composition}.
\end{proof}

\begin{proposition}\label{prop-Laplacian-generating-controlled-set}
Let $X$ be a coarse space of bounded geometry and let $\mu$ be a uniformly bounded measure.
Let $E$ be a gordo set that generates the coarse structure on $X$.
Let $E'=E^{\circ3}$.
Then for every representation $(\rho,\H)$ of $\mathbb C_{\cs}[X]$, we have $\ker(\rho(\Delta^{E'}))=\H_c$, and $X$ has geometric property (T) if and only if $\Delta^{E'}$ has spectral gap.
\end{proposition}
\begin{proof}
We apply Corollary \ref{cor-laplacians-of-close-sets-behave-the-same} with $F=E^{\circ n}$ for $n\geq 2$.
By induction it follows that for all $n\geq 3$ we have $\ker(\rho(\Delta^{E^{\circ n}}))=\ker(\rho(\Delta^{E'}))$ and that $\Delta^{E^{\circ n}}$ has spectral gap if and only if $\Delta^{E'}$ has spectral gap.

If $X$ has geometric property (T), then there is a controlled set $G$ such that $\ker(\rho(\Delta^G))=\H_c$ for every representation $(\rho,\H)$ and $\Delta^G$ has spectral gap.
Since $E$ is generating, there is some $n$ such that $G\subseteq E^{\circ n}$.
Then $\Delta^G\leq_{\max}\Delta^{E^{\circ n}}$.
So $\ker(\rho(\Delta^{E'}))=\ker(\rho(\Delta^{E^{\circ n}}))=\H_c$ for every representation $(\rho,\H)$ and $\Delta^{E^{\circ n}}$ and $\Delta^{E'}$ have spectral gap.
\end{proof}

We will now give a third characterisation of geometric property (T) based more directly on the C*-algebra $C^*_{\max}(X)$.
Note that $\ker(\Phi)$ is a left ideal in $\mathbb C_{\cs}[X]$.
Let $I_c(X)=C^*_{\max}(X)\ker(\Phi)$ be the left ideal in $C^*_{\max}(X)$ generated by $\ker(\Phi)$.

\begin{proposition}\label{prop-(T)-iff-p-exists}
Let $X$ be a coarse space of bounded geometry and let $\mu$ be a uniformly bounded measure for which a gordo set exists.
Then $(X,\mu)$ has geometric property (T) if and only if there is a projection $p\in C^*_{\max}(X)$ generating $I_c(X)$.
\end{proposition}
\begin{proof}
First suppose that $(X,\mu)$ has geometric property (T).
Let $(\rho,\H)$ be a faithful unital representation of $C^*_{\max}(X)$.
By Proposition \ref{prop-property-(T)-iff-laplacian-spectral-gap} there is a controlled set $E$ such that $\ker(\rho(\Delta^E))=\H_c$ and $\Delta^E$ has spectral gap.
By functional calculus there is a projection $p\in I_c(X)$ such that $pv=0$ for all $v\in\H_c$ and $pv=v$ for all $v\in\H_c^\perp$.
Now for any $t\in I_c(X)$ we have $tv=0$ for all $v\in\H_c$, hence $t=tp$.
So $p$ generates $I_c(X)$.

Conversely, suppose $p\in C^*_{\max}(X)$ is a projection that generates $I_c(X)$.
Then we can write $p=\sum_{j=1}^Kt_jS_j$ where $K$ is an integer and $t_j\in C^*_{\max}(X)$ and $S_j\in\ker(\Phi)\subseteq\mathbb C_{\cs}[X]$.
Choose $T'_j\in\mathbb C_{\cs}[X]$ with $\norm{t_j-T'_j}\leq \frac1{2K\norm {S_j}}$.
Let $T=\sum_{j=1}^KT_jS_j$, then $T\in\mathbb C_{\cs}[X]$, and $\Phi(T)=0$, and $\norm{T-p}\leq \sum_{j=1}^K\norm{t_j-T'_j}\cdot\norm S_j\leq \frac12$.
Write $T=T_1+\ldots+T_N$ as a sum of operators in blocks and choose a controlled set $E$ and a constant $\epsilon>0$ such that they are all in $(E,\epsilon)$-blocks.
Now let $(\rho,\H)$ be any representation.
For any $v\in\ker(p)$ we have $sv=0$ for all $s\in I_c(X)$, hence $v\in\H_c$.
So $\ker(p)=\H_c$.
Now let $v\in\H_c^\perp$ be a unit vector.
Then $pv=v$ so $\norm{Tv}\geq\frac12$.
Then there is $1\leq k\leq N$ such that $\norm{T_kv}\geq\frac1{2N}$, showing that $(X,\mu)$ has geometric property (T).
\end{proof}

\section{Coarse invariance}

In this section, we will prove that geometric property (T), as defined in Definition \ref{def-geometric-T}, does not depend on the measure $\mu$ chosen, and moreover, that it is a coarse invariant.
Our method is inspired by the one employed in \cite{WilYu14}.
Let $X$ be a coarse space with bounded geometry, and let $\mu$ be a uniformly bounded measure for which a gordo set exists.
We will start by constructing a discrete space $Y$ with bounded geometry that is coarsely dense in $X$, with an appropriate measure $\nu$.
We will show a correspondence between representations of $\mathbb C_{\cs}[Y]$ and some of the representations of $\mathbb C_{\cs}[X]$.
This will allow us to see $C^*_{\max}(Y)$ as a subalgebra of $C^*_{\max}(X)$.
Then we can apply the criterion for geometric property (T) given in Proposition \ref{prop-(T)-iff-p-exists} to conclude that $(X,\mu)$ has geometric property (T) if and only if $(Y,\nu)$ has geometric property (T).
After this we will show that geometric property (T) does not depend on  the measure as long as the space is discrete.
Finally, we will conclude that if $X$ is coarse equivalent to another space $X'$, with a uniformly bounded measure $\mu'$ for which a gordo set exists, then $(X,\mu)$ has geometric property (T) if and only if $(X',\mu')$ has geometric property (T).

By Lemma \ref{lem-complete_blocking_collection}, there is a blocking collection with union $X$.
Write $X=\bigsqcup_{y\in Y}U_y$, where $(U_y)$ is a blocking collection and $y$ is a point in $U_y$ (the notation $U_y$ is not to be confused with the notation $E_x$ introduced in Definition \ref{def-defs-coarse-set}(ii)).
Note that $Y$ is coarsely dense in $X$ and that $\mu(U_y)$ is bounded from above and below, uniformly in $y$.
Let $\nu$ be the measure on $Y$ defined by $\nu(\{y\})=\mu(U_y)$.
We will first show that $(X,\mu)$ has property (T) if and only if $(Y,\nu)$ has property (T).

\begin{lemma}\label{lem-X-T-iff-Y-T}
Let $X$ be a coarse space with bounded geometry and let $\mu$ be a uniformly bounded measure on $X$ for which a gordo set exists.
Let $Y$ be a coarsely dense discrete space of bounded geometry with a measure $\nu$.
Suppose $X=\bigsqcup_yU_y$ is a disjoint union of measurable sets with $y\in U_y$ and $\nu({y})=\mu(U_y)$, such that $\bigsqcup_yU_y\times U_y$ is a controlled set.
Then $(X,\mu)$ has property (T) if and only if $(Y,\nu)$ has property (T).
\end{lemma}

We need some preparation before we can prove Lemma \ref{lem-X-T-iff-Y-T}.
For $x\in X$, let $y(x)\in Y$ be the unique point such that $x\in U_y$.
Define the linear maps $a\colon L^2X\to L^2(Y,\nu)$ and $b\colon L^2(Y,\nu)\to L^2X$ by $a\xi(y)=\frac1{\nu(\{y\})}\int_{U_y}\xi$ and $b\eta(x)=\eta(y(x))$.
Note that these are adjoint to each other and that $ab$ is the identity on $L^2(Y,\nu)$.
Let $A=ba$.
This is a projection in $\mathbb C_{\cs}[X]$.
The maps $\alpha\colon A\mathbb C_{\cs}[X]A\to \mathbb C_{\cs}[Y],\;\alpha(T)=aTb$ and $\beta\colon \mathbb C_{\cs}[Y]\to A\mathbb C_{\cs}[X]A,\;\beta(S)=bSa$ are $*$-isomorphisms and inverse to each other.
Now if $\H^X$ is a representation of $\mathbb C_{\cs}[X]$, then $\H^Y=A\H^X$ is a representation of $A\mathbb C_{\cs}[x]A$ and via $\beta$ also a representation of $\mathbb C_{\cs}[Y]$.
Conversely, every representation of $\mathbb C_{\cs}[Y]$ gives rise to a representation of $\mathbb C_{\cs}[X]$.

\begin{lemma}\label{lem-correspondence-representations}
Let $(\rho^Y,\H^Y)$ be a representation of $\mathbb C_{\cs}[Y]$.
Then there is a representation $(\rho^X,\H^X)$ such that $A\H^X=\H^Y$, and such that $\rho^Y$ equals the restriction of $\rho^X\circ\beta$ to $A\H^X$.
Here $\H^X$ is a completion of a quotient of the algebraic tensor product $\mathbb C_{\cs}[X]\odot\H^Y$ that we will denote by $\mathbb C_{\cs}[X]\tensor\H^Y$.
\end{lemma}
\begin{proof}
First consider the algebraic tensor product $\mathbb C_{\cs}[X]\odot\H^Y$.
We equip this with the sesquilinear form defined on simple tensors by
\[\langle T\odot v,T'\odot v'\rangle=\langle\alpha(A(T')^*TA)v,v'\rangle\]
and extended to $(\mathbb C_{\cs}[X]\odot\H^Y)^2$ by sesquilinearity.
It is easy to see that this is well defined and defines a conjugate symmetric sesquilinear form.
We will now show that it is positive semi-definite.

For any operator $T\in\mathbb C_{\cs}[X]$ there is an integer $M$, such that for every $y\in Y$ there are at most $M$ different $y'\in Y$ such that $\pi_{U_y}Ti_{U_{y'}}\neq0$ or $\pi_{U_{y'}}Ti_{U_y}\neq0$.
By a graph colouring argument we may write $T=T_1+\ldots+T_N$, such that for every $1\leq k\leq N$, and every $y\in Y$, there is at most one $y'$ such that $\pi_{U_y}T_ki_{U_{y'}}\neq0$ or $\pi_{U_{y'}}T_ki_{U_y}\neq0$.

It follows that each element of $\mathbb C_{\cs}[X]\odot\H^Y$ is of the form $\sum_{k=1}^NT_k\odot v_k$ with $T_k\in\mathbb C_{\cs}[X],v_k\in\H^Y$ and the $T_k$ satisfy that for every $y$ there is at most one $y'$ such that $\pi_{U_y}Ti_{U_{y'}}\neq0$ or $\pi_{U_{y'}}Ti_{U_y}\neq 0$.
Let this $y'$ be $f_k(y)$ if it exists, otherwise put $f_k(y)=y$.
Then each $f_k$ is an involutive function (i.e. $f_k\circ f_k=\text{id}_Y$) and $T_k$ sends $L^2U_y$ to $L^2U_{f_k(y)}$.
Define
\[Q=
\begin{pmatrix}
0&T_1&\dots&T_N\\
0&0&\dots&0\\
\vdots&\vdots&\ddots&\vdots\\
0&0&\dots&0
\end{pmatrix}
\in B(L^2X)\tensor M_{N+1}(\mathbb C)=B(L^2(X\times\{0,1,\ldots,N\})).\]
For $y\in Y$ define $V_y=U_y\times\{0\}\sqcup\left(\bigsqcup_{k=1}^NU_{f_k(y)}\times\{k\}\right)\subseteq X\times\{0,1,\ldots,N\}$.
Then $X\times\{0,1,\ldots,N\}=\bigsqcup_yV_y$ and $Q$ sends each $L^2V_y$ to $L^2V_y$.
Let
\[A\tensor I=\begin{pmatrix}A&\ldots&0\\\vdots&\ddots&\vdots\\0&\dots&A\end{pmatrix}\in B(L^2(X\times\{0,\ldots,N\})).\]
Then also $(A\tensor I)Q^*Q(A\tensor I)$ sends each $L^2V_y$ to $L^2V_y$.
It is also a positive element so we can define the square root $R=\left((A\tensor I)Q^*Q(A\tensor I)\right)^\frac12 \in (A\tensor I)B(L^2(X\times\{0,1,\ldots,N\}))(A\tensor I)$.
It sends each $L^2V_y$ to $L^2V_y$.
Denote the entries of $R$ by $R_{kl}\in AB(L^2X)A$.
Then $R_{kl}$ sends $L^2U_y$ to $L^2U_{f_k(f_l(y))}$.
So $R_{kl}$ has controlled support, and we have $R_{kl}\in A\mathbb C_{\cs}[X]A$.
We have $R^2=(A\tensor I)Q^*Q(A\tensor I)$ and looking at the entry at $(k,l)$ gives $\sum_{s=1}^NR_{ks}R_{sl}=AT_k^*T_lA$ and $R_{kl}^*=R_{lk}$ for all $1\leq k,l\leq N$.

Now we can compute
\begin{align*}
&\left\langle\sum_kT_k\odot v_k,\sum_kT_k\odot v_k\right\rangle\\
=&\sum_{k,l}\langle\alpha(AT_k^*T_lA)v_l,v_k\rangle\\
=&\sum_{k,l,s}\langle\alpha(R_{ks}R_{sl})v_l,v_k\rangle\\
=&\sum_{k,l,s}\langle\alpha(R_{sl})v_l,\alpha(R_{sk})v_k\rangle\\
=&\sum_s\left\langle\sum_k\alpha(R_{sk})v_k,\sum_k\alpha(R_{sk})v_k\right\rangle\geq0.
\end{align*}

So the sesquilinear form is semi-positive definite.
Define $\H^X$ to be the Hilbert space obtained by dividing $\mathbb C_{\cs}[X]\odot\H^Y$ by the kernel of the semi-norm $\norm v=\langle v,v\rangle^\frac12$ and completing with respect to the resulting norm.
Now we will define the representation $\rho^X\colon \mathbb C_{\cs}[X]\to B(\H^X)$.

Let $T\in\mathbb C_{\cs}[X]$.
This acts on $\mathbb C_{\cs}[X]\odot\H^Y$ by $T\cdot (S\odot v)=TS\odot v$ and extending by linearity.
To show that this defines a bounded map on $\H^X$ we need to show that there is a constant $C$ depending on $T$ such that $\norm{Tw}\leq C\cdot\norm w$ for all $w\in\mathbb C_{\cs}[X]\odot\H^Y$.
For this we can assume that $T$ is an operator in blocks, by Lemma \ref{lem-sum-of-blocks}.
Then $S=(\norm{T}^2-T^*T)^\frac12$ is an element in $\mathbb C_{\cs}[X]$ and we have for all $w\in \mathbb C_{\cs}[X]\odot\H^Y$:
\[\norm T^2\cdot\norm w^2-\norm{Tw}^2=\langle (\norm T^2-T^*T)w,w\rangle=\langle Sw,Sw\rangle\geq0.\]
So $\norm{Tw}\leq \norm T\cdot\norm w$, showing that the action of $T$ on $\mathbb C_{\cs}[X]\odot\H^Y$ defines a bounded linear map $\rho^X(T)\in B(\H^X)$.
It is now easy to check that this defines a representation of $\mathbb C_{\cs}[X]$ on $\H^X$.

Finally, note that we have an injection $\H^Y\hookrightarrow\mathbb C_{\cs}[X]\odot\H^Y$ sending $v$ to $1\odot v$.
This respects the sesquilinear form, so it induces an injection $\H^Y\hookrightarrow\H^X$.
For a simple tensor $T\odot v\in \mathbb C_{\cs}[X]\odot\H^Y$ it is easy to check that $A(T\odot v)=AT\odot v=1\odot \alpha(ATA)v$, hence $A(\mathbb C_{\cs}[X]\odot\H^Y)=\H^Y$ and $A\H^X=\H^Y$.
Then it is also easy to check that the restriction of $\rho^X\circ\beta$ to $A\H^X$ equals $\rho^Y$.
\end{proof}

\begin{lemma}\label{lem-beta-extends}
The map $\beta\colon \mathbb C_{\cs}[Y]\xrightarrow\sim A\mathbb C_{\cs}[X]A$ extends to an isometry $\beta\colon C^*_{\max}(Y)\xrightarrow\sim AC^*_{\max}(X)A$.
\end{lemma}
\begin{proof}
Let $S\in\mathbb C_{\cs}[Y]$ and let $(\rho^Y,\H^Y)$ be a representation such that $\norm S_{\max}=\norm{\rho(S)}$.
Now $\H^X=\mathbb C_{\cs}[X]\tensor\H^Y$ is a representation of $\mathbb C_{\cs}[X]$ by Lemma \ref{lem-correspondence-representations}.
Moreover, the restriction of $\rho^X(\beta(S))$ to $\H^Y=A\H^X$ equals $\rho^Y(S)$, so $\norm{\beta(S)}_{\max}\geq\norm{\rho^X(\beta(S))}\geq\norm{\rho^Y(S)}_{\max}=\norm S_{\max}$.

Conversely, let $(\rho^X,\H^X)$ be a representation of $\mathbb C_{\cs}[X]$ such that $\norm{\beta(S)}_{\max}=\norm{\rho^X(\beta(S))}$.
We have $\H^X=A\H^X+(1-A)\H^X$ where $\H^Y=A\H^X$ is a representation of $\mathbb C_{\cs}[Y]$.
The restriction of $\rho^X(\beta(S))$ to $\H^Y$ equals $\rho^Y(S)$ and the restriction of $\rho^X(\beta(S))$ to $(1-A)\H^X$ equals zero, because $\beta(S)\in A\mathbb C_{\cs}[X]A$.
So $\norm S_{\max}\geq\norm{\rho^Y(S)}=\norm{\rho^X(\beta(S))}=\norm{\beta(S)}_{\max}$.

Hence $\beta\colon \mathbb C_{\cs}[Y]\to A\mathbb C_{\cs}[X]A$ is an isometry, and it extends to an isometric embedding $\beta\colon C^*_{\max}(Y)\hookrightarrow C^*_{\max}(X)$.
Its image is the closure of $A\mathbb C_{\cs}[X]A$, this is $AC^*_{\max}(X)A$.
\end{proof}

\begin{lemma}\label{lem-beta-IsY}
We have $C^*_{\max}(X)A\ker(\Phi^X)A=I_c(X)A$ and $\beta(I_c(Y))=AI_c(X)A$.
\end{lemma}
\begin{proof}
Recall that $I_c(X)=C^*_{\max}(X)\ker(\Phi^X)$, so clearly $C^*_{\max}(X)A\ker(\Phi^X)A\subseteq I_c(X)A$.
To show the converse, it is enough to show that $\ker(\Phi^X)A\subseteq \mathbb C_{\cs}[X]A\ker(\Phi^X)A$.
Let $T\in\ker(\Phi^X)A$.
For every $y,y'\in Y$ the map $\pi_{U_y}Ti_{U_{y'}}\colon L^2U_{y'}\to L^2U_y$ is the same when composed with $i_{U_{y'}}Ai_{U_{y'}}\colon L^2U_{y'}\to L^2U_{y'}$, so it has rank at most 1.
There is an integer $N$, such that for every $y\in Y$ there are at most $N$ different $y'\in Y$ such that $\pi_{U_y}Ti_{U_{y'}}\neq 0$, because $T$ has controlled support.
Therefore $\im(T)\cap L^2U_y$ has rank at most $N$.
Write $\im(T)\cap L^2U_y=f_{1y}\mathbb C \oplus \cdots\oplus f_{Ny}\mathbb C$ where the $f_{ky}$ are pairwise perpendicular elements of $L^2X$ (some of them may be zero).
For all $1\leq k\leq N$ and $y\in Y$ let $B_{ky}\colon L^2U_y\to L^2U_y$ be a rank-one isometry sending $f_{ky}\mathbb C$ to $\mathbb 1_{U_y}\mathbb C$ (or let $B_{ky}=0$ if $f_{ky}=0$).
Then $B_{ky}^*AB_{ky}=B_{ky}^*B_{ky}$ is the projection on $f_{ky}\mathbb C$.
Then $\sum_{k=1}^NB_{ky}^*AB_{ky}$ is the projection on $\im(T)\cap L^2U_y$.
Now define $B_k=\sum_yB_{ky}$.
This is a bounded operator with controlled support so we have $B_k\in\mathbb C_{\cs}[X]$.
Moreover $T=\sum_y\sum_{k=1}^NB_{ky}^*AB_{ky}=\sum_{k=1}^NB_k^*AB_kT\in\mathbb C_{\cs}[X]A\ker(\Phi^X)A$.

For the second claim note that for any $S\in \mathbb C_{\cs}[Y]$ we have $\Phi^Y(S)=0$ if and only if $\Phi^X(\beta(S))=0$.
Thus $\beta(I_c(Y))=\beta(C^*_{\max}(Y)\ker(\Phi^Y))=AC^*_{\max}(X)A\ker(\Phi^X)A=AI_c(X)A$.
\end{proof}

Now we can prove Lemma \ref{lem-X-T-iff-Y-T}.
\begin{proof}[Proof of Lemma \ref{lem-X-T-iff-Y-T}]
We will use the characterisation of geometric property (T) given in Proposition \ref{prop-(T)-iff-p-exists}.

Suppose that $(X,\mu)$ has geometric property (T).
Let $p\in C^*_{\max}(X)$ be a projection generating the left ideal $I_c(X)$.
Since $1-A\in I_c(X)$, it is a multiple of $p$, and from this it follows that $ApA$ is again a projection.
Then we see that $ApA\in AI_c(X)A$ and for all $t\in AI_c(X)A$ we have $tApA=tA=t$, so $ApA$ generates $AI_c(X)A$ as a left $AC^*_{\max}(X)A$-module.
So, by Lemma \ref{lem-beta-IsY} we know that $I_c(Y)$ is generated by a projection, hence $(Y,\nu)$ has geometric property (T).

Conversely, suppose that $(Y,\nu)$ has geometric property (T).
Let $q\in C^*_{\max}(Y)$ be a projection generating the left ideal $I_c(Y)$.
By Lemma \ref{lem-beta-IsY} we have $\beta(q)\in AI_c(X)A$.
In particular, $\beta(q)(1-A)=0$, and from this it follows that $\beta(q)+1-A$ is again a projection.
We have
\begin{align*}
&C^*_{\max}(X)(\beta(q)+1-A)\\
=&C^*_{\max}(X)\beta(q)+C^*_{\max}(X)(1-A)\\
=&C^*_{\max}(X)AI_c(X)A+C^*_{\max}(X)(1-A)\\
=&I_c(X)A+C^*_{\max}(X)(1-A)\\
=&I_c(X)
\end{align*}
by Lemma \ref{lem-beta-IsY}, so $(X,\mu)$ has geometric property (T).
\end{proof}

Now we show that geometric property (T) is independent of the measure for discrete spaces of bounded geometry.

\begin{lemma}\label{lem-different-measures-on-discrete-space}
Let $Y$ be a discrete space of bounded geometry and let $\nu$ be a measure on $Y$ such that all points are measurable and have measure uniformly bounded from above and below.
Denote the counting measure on $Y$ by $\kappa$.
Then $(Y,\nu)$ has property (T) if and only if $(Y,\kappa)$ has property (T).
\end{lemma}
\begin{proof}
Define $N\in L^\infty(Y)$ by $N(y)=\nu({y})^\frac12$.
We will consider this as a multiplication operator in $\mathbb C_{\cs}[Y]$.
Consider the linear map $a\colon L^2(Y,\nu)\to L^2(Y,\kappa)$ given by $a(\eta)=N\eta$.
This is a unitary map and it induces an isomorphism $\mathbb C_{\cs}[Y,\kappa]\to \mathbb C_{\cs}[Y,\nu]$ sending $T$ to $a^*Ta$.
For any $y\in Y$ and large enough $V\subseteq Y$, we have 
\[\Phi_\nu(a^*Ta)(y)=a^*Ta(\mathbb1_V)(y)=N^{-1}TN(\mathbb1_V)(y)=\Phi_\kappa(N^{-1}TN)(y),\]
so $\Phi_\nu(a^*Ta)=\Phi_\kappa(N^{-1}TN)$.

We now identify the algebras $\mathbb C_{\cs}[Y,\kappa]\cong\mathbb C_{\cs}[Y,\nu]\cong \mathbb C_{\cs}[Y]$, just remembering that $\Phi_\nu(T)=\Phi_\kappa(N^{-1}TN)$.
We get $C^*_{\max}(Y,\nu)=C^*_{\max}(Y,\kappa)$.
Moreover we have $\ker(\Phi_\nu)=N\ker(\Phi_\kappa)N^{-1}$ and therefore $I_c(Y,\nu)=NI_c(Y,\kappa)N^{-1}$.
We see that $I_c(Y,\nu)$ is generated by a projection if and only if $I_c(Y,\kappa)$ is generated by a projection.
By Proposition \ref{prop-(T)-iff-p-exists} we are done.
\end{proof}

Finally we prove that property (T) does not depend on the chosen measure and is also a coarse invariant.
\coarseinvariance
\begin{proof}
Construct $Y$ and $Y'$ as before.
Then $Y,X,X',Y'$ are all coarsely equivalent.
Since $Y$ and $Y'$ are coarsely equivalent, by a well-known structural result (see e.g. \cite[Lemma 4.1]{WilYu14}) there are coarse dense $Z\subseteq Y$ and $Z'\subseteq Y'$ and a bijective coarse equivalence $f\colon Z\to Z'$.
Taking the counting measure on $Z$ and $Z'$, it is easy to see that $f$ induces an isomorphism between the algebras $\mathbb C_{\cs}[Z]$ and $\mathbb C_{\cs}[Z']$ commuting with $\Phi^Z$ and $\Phi^{Z'}$, so that $Z$ has property (T) if and only if $Z'$ has property (T).

We can write $Y$ as a bounded disjoint union $\bigcup_{z\in Z}V_z$ such that $z\in V_z$.
Define the measure $\lambda$ on $Z$ by $\lambda({z})=\nu(V_z)$.
Similarly define a measure $\lambda'$ on $Z'$.
Now by Lemmas \ref{lem-X-T-iff-Y-T} and \ref{lem-different-measures-on-discrete-space} the following are equivalent:
\begin{itemize}
    \item The space $(X,\mu)$ has property (T).
    \item The space $(Y,\nu)$ has property (T).
    \item The space $(Z,\lambda)$ has property (T).
    \item The space $Z$ has property (T) with the counting measure.
    \item The space $Z'$ has property (T) with the counting measure
    \item The space $(Z',\lambda')$ has property (T).
    \item The space $(Y',\nu')$ has property (T).
    \item The space $(X',\mu')$ has property (T).
\end{itemize}
\end{proof}

\section{Amenability}

In this section, we will consider which connected coarse spaces $X$ have geometric property (T).
Recall that a coarse space $X$ is connected if $\{(x,y)\}$ is a controlled set for all $x,y\in X$.
As in \cite{WilYu14}, it turns out that an unbounded connected coarse space $X$ has geometric property (T) precisely if it does not satisfy a version of amenability.
First, we give two equivalent definitions of amenability for connected coarse spaces.
It is a generalisation of amenability for discrete metric spaces of bounded geometry, as defined e.g. in \cite{WilYu14}.

\begin{proposition}\label{prop-amenable-coarse-space}
Let $X$ be a connected coarse space and $\mu$ a uniformly bounded measure for which a gordo set exists.
The following are equivalent: 
\begin{enumerate}[(i)]
    \item There is a positive unital linear map $\phi\colon L^\infty X\to \mathbb C$ such that if $(A_i)$ is a blocking collection and $f\in L^\infty X$ satisfies $\int_{A_i}f=0$ for all $i$, and $f$ is zero outside of the union $\bigsqcup_iA_i$, then $\phi(f)=0$.
    \item For each controlled measurable set $E\subseteq X\times X$ and each $\epsilon>0$ there is a non-empty bounded measurable $U\subseteq X$ such that $\mu(E_U)\leq (1+\epsilon)\mu(U)$.
\end{enumerate}
In this case, we say that $X$ is \emph{amenable}.
\end{proposition}
\begin{proof}
First assume $(ii)$.
Define the directed set
\[\Lambda=\{(E,\epsilon)\mid E\subseteq X\times X\text{ symmetric, measurable, controlled and containing the diagonal, }\epsilon>0\}\] with $(E,\epsilon)\geq (E',\epsilon')$ if $E\supseteq E'$ and $\epsilon\leq \epsilon'$.
For each $\lambda=(E,\epsilon)\in\Lambda$, choose a measurable set $U_\lambda\subseteq X$ such that $\mu(E_{U_\lambda})\leq (1+\epsilon)\mu(U_\lambda)$.
Define $\phi_\lambda\in (L^\infty X)^*$ by $\phi_\lambda(f)=\frac1{\mu(U_\lambda)}\int_{U_\lambda}f$.
Then the $\phi_\lambda$ form a net in the unit ball of $(L^\infty X)^*$.
By the Banach-Alaoglu theorem there is a limit point $\phi$.

We will now show that $\phi$ satisfies the condition.
So let $(A_i)$ be a blocking collection and $f\in L^\infty X$ such that $\int_{A_i}f=0$ for all $i$ and $f$ is zero outside $\bigsqcup_iA_i$.
Let $E$ be the symmetric measurable controlled set $\sqcup_iA_i\times A_i\cup \{(x,x)\mid x\in X\}$.
Let $\lambda=(E',\epsilon)$ with $E'\supseteq E$.
Since $E_{U_\lambda}$ is a union of some of the $A_i$ and points outside of $\bigsqcup_iA_i$, we have $\int_{E_{U_\lambda}}f=0$.
Then $\phi_\lambda(f)=\frac1{\mu(U_\lambda)}\int_{U_\lambda}f=-\frac1{\mu(U_\lambda)}\int_{E_{U_\lambda}\setminus U_\lambda}f$ and
$|\phi_\lambda(f)|\leq \frac{\mu(E_{U_\lambda}\setminus U_\lambda)}{\mu(U_\lambda)}\cdot\norm f_\infty\leq\epsilon\norm f_\infty$.
This shows that $\phi_\lambda(f)\to0$ as $\lambda\to\infty$, hence $\phi(f)=0$.

Conversely, assume $(i)$.
Let $\phi\colon L^\infty X\to\mathbb C$ be a positive unital linear map satisfying the condition.
Let $E$ be a gordo set and assume for contradiction that there is $\epsilon>0$ such that for all bounded measurable $U\subseteq X$ we have $\mu(E_U)>(1+\epsilon)\mu(U)$.
Let $(A_i)$ be a blocking collection with union $X$, such that each $A_i$ has measure at least $\epsilon'>0$.
We may assume without loss of generality that $E=\{A_i\times A_j\mid (A_i\times A_j)\cap E\neq\emptyset\}$.
Because of bounded geometry, we can use a colouring argument to find involutions $\sigma_1,\ldots,\sigma_N$ on the index set $I$ such that $E_{A_i}=A_{\sigma_1(i)}\cup\cdots\cup A_{\sigma_N(i)}$ for all $i\in I$.
For $1\leq k\leq N$ let $E_k=\bigsqcup_iA_i\times A_{\sigma_k(i)}$.
Note that $\Delta^{E_k}$ can be viewed as an operator on $L^1$ and also as an operator on $L^\infty X$.
For any $f\in L^\infty X$ we have $\int_{A_i\cup A_{\sigma_k(i)}}\Delta^{E_k}f=0$ for all $i$, so by the condition on $\phi$ we have $\phi(\Delta^{E_k}f)=0$.
Equivalently, $(\Delta^{E_k})^*\phi=0$.

Let $\mathcal P(X)\subseteq L^1X$ denote the positive measurable functions on $X$ with integral 1.
We can view $L^1X$ as a subset of $(L^\infty X)^*$ using the standard pairing between $L^1X$ and $L^\infty X$.
By the Goldstine theorem, $\phi$ is in the weak closure of $\mathcal P(X)$.
Note that $(\Delta^{E_k})^*$ is a weakly continuous function on $(L^\infty X)^*$ and on $L^1X$ it is the same as $\Delta^{E_k}$.
Therefore, 0 is in the weak closure of the convex set
\[\{(\Delta^{E_k}\psi)\mid \psi\in\mathcal P(X)\}\subseteq\bigoplus_{k=1}^N L^1(X).\]
Since norm-closed convex sets are also weakly closed, it follows that 0 is also in the norm closure of this set.
So we can choose $\psi\in\mathcal P(X)$ such that all $\Delta^{E_k}\psi$ are small in the $L^1$-norm.
In particular, we can choose $\psi\in\mathcal P(X)$ such that $\sum_iM_i\leq\eta$, where $M_i=\sup_k\int_{A_i}|\Delta^{E_k}\psi|$ and we will choose the constant $\eta>0$ later.

Let $\delta=\frac1{1+\frac12\epsilon}<1$.
For $n\in\mathbb Z$ define $U_n=\bigsqcup\{A_i\mid\int_{A_i}\psi\geq \delta^n\mu(A_i)\}$.
Since $\int_X\psi=1$, the $U_n$ have bounded measure, they consist of finitely many of the $A_i$ and they are bounded.
Moreover, $U_n=\emptyset$ for $n$ small enough.
By assumption, we have $\mu(E_{U_n})\geq (1+\epsilon)\mu(U_n)$ for each $n$.
For each $i$, define $n_i$ as the minimal number such that $A_i\subseteq E_{U_{n_i}}$ and $N_i$ as the maximal number such that $A_i\not\subseteq U_{N_i}$ (these may be $\infty$).
Then $\int_{A_i}\psi<\delta^{N_i}\mu(A_i)$, and there is $k$ such that $\int_{A_{\sigma_k(i)}}\psi\geq \delta^{n_i}\mu(A_{\sigma_k(i)})$.
Now
\[M_i\geq \int_{A_i}|\Delta^{E_k}\psi|=\int_{A_i}\left|\mu(A_{\sigma_k(i)})\psi-\int_{A_{\sigma_k(i)}}\psi\right|\geq
\epsilon'(\delta^{n_i}-\delta^{N_i}).\]
Define
\[D_n=E_{U_n}\setminus U_{n+1}=\bigsqcup\{A_i\mid n_i\leq n\leq N_i-1\}.\]
Then
\[\sum_n\delta^n\mu(D_n)=\sum_i\mu(A_i)\sum_{n=n_i}^{N_i-1}\delta^n=\frac1{1-\delta}\sum_i\mu(A_i)(\delta^{n_i}-\delta^{N_i})\leq \frac1{(1-\delta)\epsilon'}\sum_iM_i\leq \frac\eta{(1-\delta)\epsilon'}.\]
For all $n$ we have
\[(1+\epsilon)\mu(U_n)\leq \mu(E_{U_n})\leq \mu(U_{n+1})+\mu(D_n).\]
Multiplying this inequality by $\delta^n$ on both sides and adding it for all $n$ we get
\[(1+\epsilon)\sum_n\delta^n\mu(U_n)\leq \frac1\delta\sum_n\delta^n\mu(U_n)+\sum_n\delta^n\mu(D_n),\]
so
\[(1+\epsilon-\frac1\delta)\sum_n\delta^n\mu(U_n)\leq \frac\eta{(1-\delta)\epsilon'}\]
and
\[\sum_n\delta^n\mu(U_n)\leq \frac{2\eta}{(1-\delta)\epsilon\epsilon'}.\]
On the other hand, we have
\[1=\int\psi=\sum_n\int_{U_{n+1}\setminus U_n}\psi\leq \sum_n\delta^n(\mu(U_{n+1})-\mu(U_n))=(\frac1\delta-1)\sum_n\delta^n\mu(U_n).\]
Combining these inequalities gives $1\leq \frac{2\eta}{\delta\epsilon\epsilon'}$.
Picking $\eta=\frac{\delta\epsilon\epsilon'}4$ gives the desired contradiction.
\end{proof}

The following proposition is a generalisation of \cite[Proposition 6.1]{WilYu14}.
It gives several characterisations of amenability of a coarse space $X$ in terms of its algebra.

\begin{proposition}\label{prop-amenable-coarse-space-2}
Let $X$ be a connected coarse space and $\mu$ a uniformly bounded measure for which a gordo set exists.
The following are equivalent:
\begin{enumerate}[(i)]
    \item The space $X$ is amenable.
    \item There is a net $\xi_\lambda$ of unit vectors in $L^2X$ satisfying $\norm{T\xi_\lambda-M_{\Phi(T)}\xi_\lambda}\to0$ for all $T\in\mathbb C_{\cs}[X]$ for which $\Phi(T)$ is bounded.
    \item For all $T\in\mathbb C_{\cs}[X]$ with $\Phi(T)=0$, we have $0\in\sigma(T)$.
    \item For all gordo sets $E$, we have $0\in\sigma_{\max}(\Delta^E)$.
    \item There is a unital representation $\H$ of $\mathbb C_{\cs}[X]$ that has a non-zero constant vector.
    \item There is a positive unital linear map $\phi\colon L^\infty X\to\mathbb C$ satisfying $\phi(Tf)=\phi(\overline{\Phi(T^*)}\cdot f)$ for all $T\in \mathbb C_{\cs}[X]$ for which $\phi(T^*)$ is bounded.
\end{enumerate}
\end{proposition}
\begin{proof}
First, assume that $X$ is amenable.
We will use criterion $(ii)$ from Proposition \ref{prop-amenable-coarse-space}.
As in the proof of Proposition \ref{prop-amenable-coarse-space}, define the directed set
\[\Lambda=\{(E,\epsilon)\mid E\subseteq X\times X\text{ gordo},\epsilon>0\}\] with $(E,\epsilon)\geq (E',\epsilon')$ if $E\supseteq E'$ and $\epsilon\leq \epsilon'$.
For each $\lambda=(E,\epsilon)\in\Lambda$, choose a bounded measurable subset $U_\lambda\subseteq X$ such that $\mu((E\circ E)_{U_\lambda})\leq (1+\epsilon)\mu(U_\lambda)$.
Define $\xi_\lambda=\mu(U_\lambda)^{-\frac12}\mathbb 1_{U_\lambda}$.
This is a unit vector in $L^2X$.

Let $T\in\mathbb C_{\cs}[X]$ such that $\Phi(T)$ is bounded.
We will show that $\norm{T\xi_\lambda-M_{\Phi(T)}\xi_\lambda}\to0$.
We can assume without loss of generality that $\Phi(T)=0$.
Let $T$ be supported on $E$, let $E'\supseteq E$ and $\epsilon>0$.
Then for $\lambda=(E',\epsilon)$, the function $\xi_\lambda$ is constant on $E_{X\setminus U_{\lambda}}$, so $T\xi_\lambda$ is supported on $E_{U_\lambda}$.
Since $T\mathbb 1_{(E\circ E)_{U_\lambda}}$ is zero on $E_{U_\lambda}$, we know that $T\xi_\lambda$ is the restriction to $E_{U_\lambda}$ of $-T\eta_\lambda$, where $\eta_\lambda=\mu({U_\lambda})^{-\frac12}\mathbb 1_{(E\circ E)_{U_\lambda}\setminus {U_\lambda}}$.
Now
\[\norm{T\xi_\lambda}\leq\norm{T\eta_\lambda}\leq\norm T\cdot\norm{\eta_\lambda}=\norm T\mu({U_\lambda})^{-\frac12}\mu((E\circ E)_{U_\lambda}\setminus {U_\lambda})^\frac12\leq \norm T\cdot \epsilon^\frac12,\]
so indeed $T\xi_\lambda\to0$.

The implications $(ii)\implies(iii)\implies(iv)$ are trivial.

Suppose $(iv)$ is true.
Let $\Lambda$ be the set
\[\{(E,\epsilon)\mid E\subseteq X\times X\text{ measurable controlled},\epsilon>0\}.\]
We equip this with the partial order $(E,\epsilon)\geq(E',\epsilon')$ if $E\supseteq E'$ and $\epsilon\leq\epsilon'$.
Then $\Lambda$ is a directed set.
Let $\mathcal U$ be an ultrafilter on $\mathcal P(\Lambda)$ containing $\{\lambda'\geq\lambda\}$ for each $\lambda\in\Lambda$.
By assumption, for every $\lambda=(E,\epsilon)\in\Lambda$ there is a representation $\H_\lambda$ of $\mathbb C_{\cs}[X]$ and a unit vector $v_\lambda\in\H_\lambda$ satisfying $\norm{\Delta^Ev_{\lambda}}\leq\epsilon$.
Note that for $\lambda'\geq\lambda$ we also have $\norm{\Delta^Ev_{\lambda'}}\leq\epsilon$.
Now let $\H=\prod_{\lambda\in\Lambda}\H_\lambda/\mathcal U$ be the ultraproduct of the $\H_\lambda$ and let $v=\lim_\mathcal Uv_\lambda\in\H$.
Then $v$ is a unit vector and for each measurable controlled set $E$ we have $\Delta^Ev=\lim_\mathcal U\Delta^Ev_\lambda=0$.
By Lemma \ref{lem-intersection-kernels-laplacians}, $v$ is a constant vector, showing $(iv)\implies(v)$.

For $(v)\implies(vi)$, let $v\in\H$ be a constant vector of norm 1 and define $\phi\colon L^\infty X\to\mathbb C$ by $\phi(f)=\langle M_fv,v\rangle$.
Then $\phi$ is a positive unital linear map and for all $f\in L^\infty X$ and $T\in\mathbb C_{\cs}[X]$ such that $\Phi(T^*)$ is bounded, we have
\[\phi(Tf)=\langle M_{\Phi(TM_f)}v,v\rangle=\langle TM_fv,v\rangle=\langle M_fv,M_{\Phi(T^*)}v\rangle=\langle M_{\overline{\Phi(T^*)}\cdot f}v,v\rangle=\phi(\overline{\Phi(T^*)}\cdot f).\]

Finally, for $(vi)\implies(i)$, let $\phi\colon L^\infty X\to\mathbb C$ be a positive unital linear map satisfying $\phi(Tf)=\phi(\overline{\Phi(T^*)}\cdot f)$ for all $T\in\mathbb C_{\cs}[X]$ such that $\Phi(T^*)$ is bounded.
Let $(A_i)$ be a blocking collection and $E=\bigsqcup_iA_i\times A_i$.
Let $f\in L^\infty X$ be a function such that $\int_{A_i}f=0$ for all $i$, and $f$ is zero outside of $\bigsqcup_iA_i$.
Then $\Delta^Ef=f$ and $\phi(\Delta^E)=0$, so $\phi(f)=\phi(\Delta^Ef)=0$.
Therefore, $X$ is amenable (using part $(i)$ of Proposition \ref{prop-amenable-coarse-space}).

\end{proof}

\begin{remark}
If $E'$ is as in proposition \ref{prop-Laplacian-generating-controlled-set}, the properties in Theorem \ref{prop-amenable-coarse-space-2} are also equivalent to $0\in\sigma_{\max}(\Delta^{E'})$.
\end{remark}

Now we can give the generalisation of \cite[Corollary 6.1]{WilYu14}.

\amenability
\begin{proof}
If $X$ is bounded, we have already shown that it has geometric property (T), and it is amenable because $\mathbb1_X\in L^2X$ is a constant vector in the standard representation.

Suppose that $X$ is unbounded and not amenable.
Then there is a gordo set $E$ such that $0\not\in\sigma_{\max}(\Delta^E)$.
Then for all representations $(\rho,\H)$ of $\mathbb C_{\cs}[X]$ we have $\ker(\rho(\Delta^E))=\H_c=\emptyset$ and $\Delta^E$ has spectral gap.
By proposition \ref{prop-property-(T)-iff-laplacian-spectral-gap} the space $X$ has geometric property (T).

Suppose that $X$ is unbounded, has geometric property (T) and is amenable.
By Proposition \ref{prop-property-(T)-iff-laplacian-spectral-gap} there is a measurable controlled set $E$ such that for all representations $(\rho,\H)$ of $\mathbb C_{\cs}[X]$, we have $\H_c=\ker(\rho(\Delta^E))$, and $\Delta^E$ has spectral gap.
Since $X$ is amenable, we have $0\in\sigma(\Delta^E)$, but then 0 is an isolated point in the spectrum of $\Delta^E\in B(L^2X)$.
So there is a non-zero $\xi\in\ker(\Delta^E)$, which is then a constant vector.
Since $\xi\in L^2X$ and $X$ is unbounded, $\xi$ can not be constant as a function.
So there are measurable sets $A,B\subseteq X$ of positive measure such that the convex hulls of $\xi(A)$ and $\xi(B)$ are disjoint.
If we intersect $A$ with all of the members of a blocking collection, one of the intersections must have positive measure, since blocking collections are countable.
Therefore, we can assume that $A$ is bounded, and similarly, that $B$ is bounded.
Since $X$ is connected, $A\cup B$ is also bounded.
Let $T\in\mathbb C_{\cs}[X]$ be the operator defined by
\[T\eta(x)=\begin{cases}\frac1{\mu(B)}\int_B\eta-\frac1{\mu(A)}\int_A\eta&\text{if }x\in A\\0&\text{otherwise.}\end{cases}\]
Then $\Phi(T)=0$ but $T\xi\neq 0$, hence $\xi$ is not a constant vector and we have a contradiction.
\end{proof}

\begin{corollary}
Amenability for connected coarse spaces of bounded geometry does not depend on the measure chosen and is a coarse invariant.
\end{corollary}
\begin{proof}
This follows from Theorem \ref{thm-coarse-invariance} and Theorem \ref{thm-connected-space-amenable-or-(T)}.
It is also straightforward to show it directly, using criterion $(ii)$ of Proposition \ref{prop-amenable-coarse-space}.
\end{proof}

\section{Manifolds}\label{sec-manifolds}

Let $(M,g)$ be an $n$-dimensional complete Riemannian manifold, not necessarily connected.
Consider the geodesic distance $d$ and volume $\mu$.
The metric $d$ defines a coarse structure on $M$, where $E\subseteq M\times M$ is controlled if $\sup_{(x,y)\in E}d(x,y)<\infty$.
For any $R>0$ let $E_R$ be the controlled set $\{(x,y)\in M\times M\mid d(x,y)\leq R\}$.
Let $B(x,R)=(E_R)_x$ denote the ball of radius $R$ around a point $x\in M$.
In \cite{Gri99}, it is defined that $M$ has \emph{bounded geometry} if the injectivity radius is positive, and the Ricci curvature is bounded below (by a possibly negative constant).
In this section, we assume that this is the case.

By Bishop's inequality (see \cite[Theorem 3.101.i]{GaHuLa80}) there are constants $B,B'>0$ such that $\mu(B(x,R))\leq B\exp(B'R)$ for every $x\in M$ and $R>0$.
Hence, the measure $\mu$ is uniformly bounded.
There are also $R,A>0$ such that $\mu(B(x,R))\geq A$ for all $x\in M$, by \cite[Proposition 14]{Croke80} (the Proposition is stated for compact manifolds, but the proof works just as well in the non-compact case if the injectivity radius is positive).
In particular, the controlled set $E_R$ is gordo, so we can conclude that the manifold $M$ has bounded geometry as a coarse space.

Let $\Delta_M$ be the (positive) Laplacian operator on $M$ and let $\Delta_H=1-\exp(-\Delta_M)$.
Since the spectrum of $\Delta_M$ is contained in $[0,\infty]$, the spectrum of $\Delta_H$ is contained in $[0,1]$.
In particular, $\Delta_H$ is a bounded operator.
We can write
\[\Delta_Hf(x)=f(x)-\int_Mp(x,y)f(y)d\mu(y)\]
where $p(x,y)$ is the heat kernel on $M$.
We need some estimates on this function.

\begin{lemma}\label{lem-estimates-heat-kernel}
\begin{enumerate}[(i)]
    \item There are constants $c,R>0$ such that $p(x,y)\geq c$ whenever $d(x,y)\leq R$.
    \item The manifold is \emph{stochastically complete}: for every $x\in M$ we have $\int_Mp(x,y)d\mu(y)=1$.
    \item For every $\epsilon>0$ there is an $R$ such that $\int_{B(x,R)}p(x,y)d\mu(y)\geq 1-\epsilon$ for all $x\in M$.
\end{enumerate}
\end{lemma}
\begin{proof}
\begin{enumerate}[(i)]
    \item This follows from \cite[Equation 7.43]{Gri99}.
    \item This follows from \cite[Theorem 3.4]{Gri99} and the bound on $\mu(B(x,R))$ given above.
    \item The conditions of \cite[Equation 6.40]{Gri99} are met, so there are constants $C,D>0$ such that
    \[p(x,y)\leq C\exp\left(\frac{-d(x,y)^2}D\right)\]
    for all $x,y\in M$.
    Now for any $x\in M$ and $R>0$ we have
    \begin{align*}
        \int_{M\setminus B(x,R)}p(x,y)d\mu(y)&\leq C\int_{M\setminus B(x,R)}\exp\left(\frac{-d(x,y)^2}D\right)d\mu(y)\\
        &\leq C\int_R^\infty\frac{2r}D\exp\left(\frac{-r^2}D\right)V_r(x)dr\\
        &\leq \frac{2BC}D\int_R^\infty r^{n+1}\exp\left(\frac{-r^2}D\right)dr.
    \end{align*}
    As this last integral converges, we know that for every $\epsilon>0$ there is an $R>0$ such that $\int_{M\setminus B(x,R)}p(x,y)d\mu(y)\leq\epsilon$.
    Then by part (ii) we get $\int_{B_{x,R}}p(x,y)d\mu(y)\geq 1-\epsilon$ for all $x\in M$.
\end{enumerate}
\end{proof}

\begin{lemma}\label{lem-Delta-H-in-algebra}
We have $\Delta_H\in C^*_{\max}(M)$.
\end{lemma}
\begin{proof}
For any $R>0$ let $\Delta_{HR}\in C_{\cs}[M]$ be the operator defined by
\[\Delta_{HR}f(x)=\int_{B(x,R)}p(x,y)(f(x)-f(y))d\mu(y).\]
Then for $R'>R$ we have
\[(\Delta_{HR'}-\Delta_{HR})f(x)=\int_{B(x,R')\setminus B(x,R)}p(x,y)(f(x)-f(y))d\mu(y).\]
For large enough $R$ we have $\int_{M\setminus B(x,R)}p(x,y)d\mu(y)\leq \epsilon$.
Then we have
\[\langle(\Delta_{HR'}-\Delta_{HR})\xi,\xi\rangle=\frac12\int_{\{(x,y)\in M^2\mid R<d(x,y)\leq R'\}}p(x,y)|\xi(x)-\xi(y)|^2d(\mu\times\mu)(x,y)\leq2\epsilon\norm\xi^2.\]
So
\[0\leq_{\max}\Delta_{HR'}-\Delta_{HR}\leq_{\max}2\epsilon.\]
Therefore the $\Delta_{HR}$ form a Cauchy sequence in $\mathbb C_{\cs}[M]$ with the maximal norm, and they converge to $\Delta_H$.
\end{proof}

For any $R>0$ let $\Delta_R=\Delta^{E_R}$.
We will use the estimates of the heat kernel from Lemma \ref{lem-estimates-heat-kernel} to prove estimates on $\Delta_H$ in terms of the $\Delta_R$.

\begin{lemma}\label{lem-estimates-Delta-H}
We have the following estimates on $\Delta_H$:
\begin{enumerate}[(i)]
    \item There are constants $c,R>0$ such that $\Delta_H\geq_{\max} c\Delta_R$.
    \item For every $\epsilon>0$ there are constants $d,R>0$ such that $\Delta_H\leq_{\max} \epsilon+d\Delta_R$.
\end{enumerate}
\end{lemma}
\begin{proof}
\begin{enumerate}[(i)]
    \item Take the same constants $c,R$ as in Lemma \ref{lem-estimates-heat-kernel} part (i).
    Then
    \[(\Delta_H-c\Delta_R)f(x)=\int_M\alpha(x,y)(f(x)-f(y))d\mu(y),\]
    where $\alpha(x,y)=p(x,y)$ if $d(x,y)>R$ and $\alpha(x,y)=p(x,y)-c$ if $d(x,y)\leq R$.
    Since $\alpha\geq 0$ everywhere we have $\Delta_H\geq_{\max} c\Delta_R$.
    \item Let $\Delta_{HR}$ be as in the proof of Lemma \ref{lem-Delta-H-in-algebra}.
    There is a large enough $R$ such that $\Delta_H\leq_{\max}\epsilon+\Delta_{HR}$.
    The function $p(x,y)$ is bounded from above by some constant $d$, and then we get $\Delta_{HR}\leq_{\max}d\Delta_R$.
\end{enumerate}
\end{proof}

For any representation $\H$ of $C^*_{\max}(M)$, we can consider the unbounded symmetric operator $\Delta_M=-\log(1-\Delta_H)$ on $\H$.
We say that $\Delta_M$ has spectral gap if there is a constant $\gamma>0$ such that the spectrum of this operator is contained in $\{0\}\cup[\gamma,\infty)$ for each representation $\H$.

\manifolds
\begin{proof}
For any $R>0$, the controlled set $E_R$ generates the coarse structure on $M$ and is gordo.
So we can apply Proposition \ref{prop-Laplacian-generating-controlled-set} to see that $M$ has geometric property (T) if and only if $\Delta_R$ has spectral gap.
Clearly, $\Delta_H$ has spectral gap if and only if $\Delta_M$ has spectral gap.
For the remainder of the proof, we use the estimates of Lemma \ref{lem-estimates-Delta-H}.

Suppose that $M$ has geometric property (T).
By Lemma \ref{lem-estimates-heat-kernel} part (i), there are constants $c,R>0$ such that $\Delta_H\geq c\Delta_R$.
Now $\Delta_R$ has spectral gap, and so does $\Delta_H$, and also $\Delta_M$.

Suppose that $\Delta_H$ has spectral gap.
There is $\epsilon>0$ such that $\sigma_{\max}(\Delta_H)\subseteq \{0\}\cup [2\epsilon,\infty)$.
By Lemma \ref{lem-estimates-heat-kernel} part (ii) there are $b,R>0$ such that $\Delta_H\leq \epsilon+b\Delta_R$.
Then $\Delta_R$ has spectral gap, so $M$ has geometric property (T).
\end{proof}

\section{Warped systems}\label{sec-warped-systems}

Let $\Gamma$ be a finitely generated group with finite generating set $S$.
Let $(M,g)$ be a compact Riemannian manifold.
Let $d$ be the distance function on $M$ and $\mu$ the measure on $M$ defined by the metric $g$.
Let $\alpha\colon \Gamma\curvearrowright M$ be an action by diffeomorphisms.
Let $\mathcal CM=M\times\{1,2,\ldots\}$ be the discrete cone.
This is equipped with the Riemannian metric $g_{\mathcal C}$ that is $t\cdot g$ on $M_t=M\times\{t\}$.
Let $d_\mathcal C$ be the corresponding distance function, given by $d_\mathcal C(x,y)=t\cdot d(x,y)$ when $x,y\in M_t=M\times\{t\}$ and $d_\mathcal C(x,y)=\infty$ if $x,y$ are in different $M_t$.
Now define the \emph{warped system} $\Warp(\Gamma\curvearrowright M)$ as the space $\mathcal CM$ equipped with the largest metric $\delta_\Gamma$ for which $\delta_\Gamma(x,y)\leq d_\mathcal C(x,y)$ and $\delta_\Gamma(x,s\cdot x)\leq 1$ for all $s\in \Gamma$ and $x,y\in\mathcal CM$.
We will consider $\Warp(\Gamma\curvearrowright M)$ as a coarse space.

In \cite{Vig19}, it is shown that that $X=\Warp(\Gamma\curvearrowright M)$ does not depend on the generating set of $\Gamma$ as a coarse space.
It is also coarsely equivalent to the coarse disjoint union $\bigsqcup_{t\in\mathbb Z}M\times\{t\}$ of subspaces of the warped cone as introduced by Roe in \cite{Roe05} (see \cite[Lemma 6.5]{Vig19}).

Let $\mu$ be the measure on $X$ defined through the Riemannian metric $g_\mathcal C$.
Since the manifold $M$ is compact, the injectivity radius is positive and the Ricci curvature is bounded below.
It follows as in the previous section that $F_R=\{(x,y)\in X\times X\mid d_\mathcal C(x,y)\leq R\}$ is gordo for any $R>0$, and that its balls are uniformly bounded.
By Lemma \ref{lem-controlled_set_warped_cone} below, it follows that $\mu$ is uniformly bounded.
In particular, $X$ has bounded geometry.

Let $\Delta_\Gamma=\sum_{s\in S}1-s \in \mathbb C[\Gamma]$.
This can also be viewed as an element in $\mathbb C_{\cs}[X]$.
In \cite{Vig19}, it is shown that a coarsely dense subset of $X$ with discrete bounded geometry, is a family of expanders if and only if there is a constant $\gamma>0$ such that $\langle\Delta_\Gamma \xi,\xi\rangle\geq \gamma\norm\xi^2$ for all $\xi\in L^2M$ with $\int_M\xi d\mu=0$.
We want to consider when these graphs have geometric property (T), or equivalently by Theorem \ref{thm-coarse-invariance}, when $X$ has geometric property (T).

\begin{definition}
Let $\Gamma$ be a group and $(M,\mu)$ a measure space.
Let $\alpha\colon \Gamma\curvearrowright M$ be a measurable action.
The action is called \emph{ergodic} if the only measurable $\Gamma$-invariant subsets of $M$ are the empty set and $M$ itself.
\end{definition}

If $\Gamma$ has property (T) and the action on $M$ is ergodic, we might expect the warped system $X$ to have geometric property (T).
At the moment, this is an open question.

\begin{question}
Suppose $\Gamma$ has property (T) and the action $\alpha\colon \Gamma\curvearrowright M$ is ergodic.
Does $X$ have geometric property (T)?
\end{question}
The measure-preserving action of $\Gamma$ on $M$ gives a map $\Gamma\to U(L^2X)$ given by $g\cdot \xi(x)=\xi(g^{-1}x)$ for $g\in\Gamma,\xi\in L^2X$.
The image is contained in $\mathbb C_{\cs}[X]$.
Then the map $\Gamma\to\mathbb C_{\cs}[X]$ gives rise to a map $\mathbb C[\Gamma]\rightarrow\mathbb C_{\cs}[X]$.
Now every representation of $\mathbb C_{\cs}[X]$ pulls back to a representation of $\mathbb C[\Gamma]$, therefore the map $\mathbb C[\Gamma]\to \mathbb C_{\cs}[X]$ is contracting for the maximal norms.
So it induces a map $C^*_{\max}(\Gamma)\to C^*_{\max}(X)$.
If $\Gamma$ has property (T), then we know that $\Delta_\Gamma$ has spectral gap in $C^*_{\max}(\Gamma)$.
So the image of $\Delta_\Gamma$ also has spectral gap in $C^*_{\max}(X)$.
If $\ker(\rho(\Delta))=\H_c$ for any representation $(\rho,\H)$ of $\mathbb C_{\cs}[X]$, then we can conclude that $X$ has property (T) by Proposition \ref{prop-property-(T)-iff-laplacian-spectral-gap}.
We can prove the following partial result.
\begin{lemma}
Suppose that $\Gamma$ has property (T) and the action $\alpha\colon\Gamma\curvearrowright M$ is ergodic.
Let $\H$ be a representation of $\mathbb C_{\cs}[X]$, and $v\in\H$ such that $\Delta_\Gamma v=0$.
For any $f\in L^\infty X$ such that $\int_{X_t}fd\mu=0$ for all $t\in\mathbb N$, we have $\langle M_fv,v\rangle=0$.
\end{lemma}
\begin{proof}
We can assume without loss of generality that $v$ is a unit vector.

The composition $\Gamma\xrightarrow\alpha\mathbb C_{\cs}[X]\xrightarrow\rho B(\H)$ gives an action of $\Gamma$ on $\H$, still denoted by $\alpha$.
For any $s\in S$, we have $\norm{(1-\alpha_s)v}^2=\langle (1-\alpha_s^*)(1-\alpha_s)v,v\rangle=\langle (2-\alpha_s-\alpha_s^*)v,v\rangle$.
Since $2-\alpha_s-\alpha_s^*\leq \Delta_\Gamma$ it follows that $(1-\alpha_s)v=0$.
So $v$ is $\Gamma$-invariant.

Now define $\phi\in (L^\infty X)^*$ by $\phi(f)=\langle M_fv,v\rangle$.
Note that $\phi$ is a mean on $L^\infty X$.
Moreover, $\phi$ is $\Gamma$-invariant:  for any $g\in \Gamma,f\in L^\infty X$ we have
\[\phi(\alpha_gf)=\langle M_{\alpha_gf}v,v\rangle=\langle \alpha_gM_f\alpha_g^*v\rangle=\langle M_fv,v\rangle=\phi(f).\]
Let $\mathcal P(X)=\{\psi\in L^1X \mid \psi\geq 0,\int_X\psi d\mu=1\}$.
By the Goldstine theorem, $\phi$ is in the weak closure of $\mathcal P(X)$.
Let $f\in L^\infty X$ such that $\int_{X_t}fd\mu=0$ for all $t\in\mathbb N$.
Let $\epsilon>0$.
Then 0 is in the weak closure of the convex set
\[\left\{((1-\alpha_s)\psi)\in\sum_{s\in S}L^1X\mid \psi\in\mathcal P(X),|\psi(f)-\phi(f)|<\epsilon\right\}.\]
Since norm-closed convex sets are also weakly closed, 0 is also in the norm closure of the set above.
Hence there is $\psi\in\mathcal P(X)$ with $|\psi(f)-\phi(f)|<\epsilon$ and $\norm{(1-\alpha_s)\psi}_1<\epsilon$ for all $s\in S$.
Now let $\xi=\psi^\frac12\in L^2X$.
This is a unit vector satisfying
\[\norm{(1-\alpha_s)\xi}^2_2=\int_X(\xi(s\cdot x)-\xi(x))^2d\mu(x)\leq \int_X|\xi(s\cdot x)^2-\xi(x)^2|=\norm{(1-\alpha_s)\psi}_1<\epsilon\]
for all $s\in S$.
Then $\norm{\Delta_\Gamma \xi}_2\leq \epsilon^\frac12|S|$.
Let $\xi_c\in L^2X$ be the locally constant function defined by $\xi_c(x)=\frac1{\mu(X_t)}\int_{X_t}\xi$ for $x\in X_t$.
There is $\gamma>0$ such that $\sigma(\Delta_\Gamma)\subseteq\{0\}\cup[\gamma,\infty)$.
Now $\xi-\xi_c$ is perpendicular to the locally constant functions, and the locally constant functions are the only $\Gamma$-invariant functions in $L^2X$, since the action of $\Gamma$ is ergodic.
From $\norm{\Delta_\Gamma(\xi-\xi_c)}_2\leq \epsilon^\frac12|S|$ it now follows that $\norm{(\xi-\xi_c)}_2\leq \gamma^{-1}\epsilon^\frac12|S|$.
Now
\begin{align*}
|\phi(f)|&\leq|\psi(f)|+\epsilon\\
&=\left|\int_X\xi(x)^2f(x)d\mu(x)\right|+\epsilon\\
&=\left|\int_X(\xi(x)-\xi_c(x))^2f(x)d\mu(x)\right|+2\left|\int_X(\xi(x)-\xi_c(x))\xi_c(x)f(x)d\mu(x)\right|+\epsilon\\
&\leq \gamma^{-2}\epsilon|S|^2\cdot\norm f_\infty+\gamma^{-1}\epsilon^\frac12|S|\cdot\norm f_\infty+\epsilon.
\end{align*}
Since this holds for any $\epsilon>0$, we conclude that $\phi(f)=0$.
\end{proof}

\begin{remark}
In the above, we only used that $\sigma(\Delta_\Gamma)\subseteq\{0\}\cup[\gamma,\infty)$ instead of the stronger condition $\sigma_{\max}(\Delta_\Gamma)\subseteq\{0\}\cup[\gamma,\infty)$.
\end{remark}
\begin{remark}
To see why the Lemma is progress to the Question, consider that we have to prove that every $v\in\ker(\rho(\Delta_\Gamma))$ is constant.
For this we can assume that $v\in \H_c^\perp$, and try to show that it is in fact 0.
Now from the Lemma it follows that in this case, we have $v\in (\rho(L^\infty X)\H_c)^\perp$, in other words, $v$ is a vector that is in some way quite far from being constant.

To see this, consider $f\in L^\infty X$ and $w\in\H_c$.
We can write $f=f_c+g$ with $f_c$ constant on each component of $X$ and $\int_{X_t}gd\mu=0$ for all $g$.
For all $t\in\mathbb C$ we have $v+tw\in\ker(\rho(\Delta_\Gamma))$, so by the Lemma, we have $\langle M_g(v+tw),v+tw\rangle = 0$.
Using some different values of $t$ gives $\langle v, M_gw\rangle=0$.
Also $\langle v,M_{f_c}w\rangle=0$ since $M_{f_c}w$ is a constant vector.
So $\langle v,M_fw\rangle=0$.
\end{remark}

The converse is also an open question.
\begin{question}
Suppose the action $\alpha$ is free and $X$ has geometric property (T).
Is it necessary that $\Gamma$ has property (T) and that the action is ergodic?
\end{question}

We have some partial results on this question.
Before we state them, we give some preliminary results.

\begin{lemma}\label{lem-controlled_set_warped_cone}
For $R>0$ and $g\in \Gamma$ define the controlled sets $E_R=\{(x,y)\in X\mid \delta_\Gamma(x,y)\leq R\},\,E_g=\{(g\cdot x,x)\mid x\in X\}$ and $F_R=\{(x,y)\in X\mid d_\mathcal C(x,y)\leq R\}$.
For any $R>0$ there is $R'>0$ such that $E_R\subseteq \bigcup_{g\in S^{\lfloor R\rfloor}}F_{R'}\circ E_g$.
\end{lemma}
\begin{proof}
For any $s\in S$ the map $M\to M,\;x\to s\cdot x$ is Lipschitz.
So there is a constant $L>0$ such that $d_\mathcal C(s\cdot x,s\cdot y)\leq L\cdot d_\mathcal C(x,y)$.
Let $R'=L^RR$.
For any $(y,x)\in E_R$, there are $x_0=x,x_1,\ldots,x_n \in X$ and $y_0,y_1,\ldots,y_n=y\in X$ and $s_1,s_2,\ldots,s_n\in S$ such that $s_ky_{k-1}=x_k$ for $k=1,2,\ldots,n$ and
\[n+d_\mathcal C(x_0,y_0)+\ldots+d_\mathcal C(x_n,y_n)\leq R.\]
It follows by induction that
\[d_\mathcal C(s_k\cdots s_2\cdot s_1x_0,y_k)\leq L^kd_\mathcal C(x_0,y_0)+\ldots+Ld_\mathcal C(x_{k-1},y_{k-1})+d_\mathcal C(x_k,y_k).\]
In particular
\[d_\mathcal C(g\cdot x,y)\leq L^nd_\mathcal C(x_0,y_0)+\ldots+Ld_\mathcal C(x_{n-1},y_{n-1})+d_\mathcal C(x_n,y_n)\leq L^RR=R'\]
for $g=s_n\cdots s_2\cdot s_1$.
Now $(g\cdot x,x)\in E_g$ and $(y,g\cdot x)\in F_{R'}$ so $(y,x)\in F_{R'}\circ E_g$.
Since $g\in S^{\lfloor R\rfloor}$ this shows that \[E_R\subseteq\bigcup_{g\in S^{\lfloor R\rfloor}}F_{R'}\circ E_g.\]
\end{proof}

Note that $I=\bigoplus_{t\in\mathbb N}B(L^2X_t)$ is a two-sided ideal in $\mathbb C_{\cs}[X]$.
\begin{lemma}
Suppose that the action $\alpha$ is free.
Then $\mathbb C_{\cs}[X]/I$ is naturally isomorphic to the algebraic cross product $\mathbb C_{\cs}[\mathcal CM]/I\rtimes \Gamma$.
\end{lemma}
\begin{proof}
The natural map $\mathbb C_{\cs}[\mathcal CM]\rtimes\Gamma\to \mathbb C_{\cs}[X]$ descends to $\mathbb C_{\cs}[\mathcal CM]/I\rtimes\Gamma\to\mathbb C_{\cs}[X]/I$.

To show that it is surjective, let $T\in \mathbb C_{\cs}[X]$.
Then $T$ is supported on $E_R$ for some integer $R>0$.
By Lemma \ref{lem-controlled_set_warped_cone}, there is $R'>0$ such that $E_R\subseteq \bigcup_{g\in S^R}F_{R'}E_g$.
We can use this decomposition to write $T$ as a finite sum $T=\sum_gA_g\alpha_g$, with $A_g\in \mathbb C_{\cs}[\mathcal CM]$ showing that the map is indeed surjective.

For injectivity, suppose that $\sum_{g\in B}A_g\alpha_g \in I$, with $B\subset\Gamma$ finite and $A_g\in\mathbb C_{\cs}[\mathcal CM]$.
Then we will show that $A_g\in I$ for all $g\in B$.
We may as well assume that $\sum_{g\in B}A_g\alpha_g=0$.
Since the action is free and by homeomorphisms, for any $g\in\Gamma$ we have $\min_{x\in M}d(g\cdot x,x)>0$.
There is some $R>0$ such that all $A_g$ are supported on $F_R$.
Now let $T=\frac {2R}{\min_{g\neq h\in B}\min_{x\in M}d(gh^{-1}\cdot x,x)}$.
For any $g\in B$ we know that $A_g\alpha_g$ is supported on $F_RE_g$.
We also have $A_g\alpha_g=-\sum_{h\in B\setminus \{g\}}A_h\alpha_h$, so $A_g\alpha_g$ is supported on $F_RE_g\cap\bigcup_{h\in B\setminus \{g\}}F_RE_h$.
Suppose that $(x,y)\in F_RE_g\cap F_RE_h$.
Then $d_\mathcal C(x,g\cdot y)\leq R$ and $d_\mathcal C(x,h\cdot y)\leq R$, so $d_\mathcal C(g\cdot y,h\cdot y)\leq 2R$.
For some $t$ we have $x,y\in X_t$ and then $d_\mathcal C(g\cdot y,h\cdot y)\geq t\cdot\min_{z\in M}d(gh^{-1}\cdot z,z)\geq\frac{2Rt}T$.
So $t\leq T$.
Hence we have $F_RE_g\cap F_RE_h\subseteq \bigcup_{t\leq T}X_t\times X_t$, and $F_RE_g\cap\bigcup_{h\in B\setminus\{g\}}F_RE_h\subseteq\bigcup_{t\leq T}X_t\times X_t$.
So $A_g\alpha_g$ is supported on $\bigcup_{t\leq T}X_t\times X_t$, and we conclude that $A_g\in\bigoplus_{t\in\mathbb N}B(L^2X_t)$.
\end{proof}

\begin{lemma}
Suppose the action $\alpha$ is free.
Then the natural map $C^*_{\max}(\Gamma)\to C^*_{\max}(X)$ is an injection.
\end{lemma}
\begin{proof}
The proof is similar to the proof of \cite[Lemma 7.1]{WilYu14}.
Let $\overline I$ denote the closure of $I$ in $C^*_{\max}(X)$.
The previous lemma gives $\mathbb C_{\cs}[X]/I \cong \mathbb C_{\cs}[\mathcal CM]/I\rtimes \Gamma$, hence $C^*_{\max}(X)/\overline I\cong C^*_{\max}(\mathcal CM)/\overline I\rtimes_{\max}\Gamma$.
Now it suffices to prove that the composition $C^*_{\max}(\Gamma)\to C^*_{\max}(X)\to C^*_{\max}(X)/\overline I\cong C^*_{\max}(\mathcal CM)/\overline I\rtimes_{\max}\Gamma$ is injective.

For any $t\in\mathbb N$ define the map $\phi_t\colon \mathbb C_{\cs}[\mathcal CM]\to \mathbb C$ by $\phi_t(T)=\frac1{\mu(X_t)}\langle\mathbb 1_{X_t},T\mathbb 1_{X_t}\rangle$.
This extends to a positive unital map $\phi_t\colon C^*_{\max}(\mathcal CM)\to\mathbb C$.
Let $\phi$ be a limit point of these maps.
Then $\phi\colon C^*_{\max}(\mathcal CM)/\overline I\to\mathbb C$ is a ucp map.
Since maximal crossed products are functorial for ucp maps, this gives a map $C^*_{\max}(\mathcal CM)/\overline I\rtimes_{\max}\Gamma\to C^*_{\max}(\Gamma)$.
Now the composition $C^*_{\max}(\Gamma)\to C^*_{\max}(\mathcal CM)/\overline I\rtimes_{\max}\Gamma\to C^*_{\max}(\Gamma)$ is the identity, so the first map is an injection.
\end{proof}

As a result, we get that for free actions, $\Delta_\Gamma$ has spectral gap in $C^*_{\max}(X)$ if and only if $\Gamma$ has property (T).
So if $\Gamma$ does not have property (T), for any $\epsilon>0$, there is some representation $(\rho,\H)$ and a unit vector $v\in\H_c^\perp$ with $\norm{\Delta_\Gamma v}<\epsilon$.
Suppose that $\norm{\Delta_Rv}<1$ for some $0<R<1$.
Then one can show that $\norm{\Delta_Rv}\to 0$ as $R\to 0$.
Let $u\colon L^2X\to L^2X$ be given by $u(\xi)(x,t)=\xi(x,2t)$.
Precomposing the representation $\rho$ by the map $\tau\colon \mathbb C_{\cs}[X]\to\mathbb C_{\cs}[X], \tau(T)=uTu^*$ multiple times gives vectors $v_0=v,v_1,\ldots$ such that $\norm{\Delta_\Gamma v_n}<\epsilon$ and $\norm{\Delta_R v_n}<\epsilon$.
In this case it follows that $X$ does not have geometric property (T).
A similar strategy would work for non-ergodic actions.
The main obstacle is, that there may be vectors $v$ for which $\norm{\Delta_Rv}=1$ for all $0<R<1$.

\section*{Acknowledgement}
I would like to thank my advisor Tim de Laat for his support and suggestions.
I also thank for Rufus Willett for helpful suggestions, and the reviewer for a very thorough and helpful report.
The author is supported by the Deutsche Forschungsgemeinschaft under Germany's Excellence Strategy - EXC 2044 - 390685587, Mathematics M\"unster: Dynamics - Geometry - Structure.


\begin{thebibliography}{99}

\bibitem{Croke80}
C. B. Croke, \emph{Some isoperimetric inequalities and eigenvalue estimates}, Annales scientifiques de l’É.N.S. \textbf{13}(1980), 419-435

\bibitem{GaHuLa80}
S. Gallot, D. Hulin and J. Lafontaine, \emph{Riemannian geometry}, Springer-Verlag, Berlin, 1980

\bibitem{Gri99}
A. Grigor'yan, \emph{Estimates of heat kernels on Riemannian manifolds}, Cambridge University Press, Cambridge, 1999

\bibitem{HigPedRoe95}
N. Higson, E. Pedersen and J. Roe, \emph{C*-algebras and controlled topology}, K-theory \textbf{11}(1997), 209-239

\bibitem{Roe03}
J. Roe, \emph{Lectures on coarse geometry}, Amer. Math. Soc., Providence, RI, 2003

\bibitem{Roe05}
J. Roe, \emph{Warped cones and property A}, Geometry and Topology \textbf{9}(2005), 163-178

\bibitem{Vig19}
F. Vigolo, \emph{Measure expanding actions, expanders and warped cones}, Trans. Amer. Math. Soc. \textbf{371}(2019), 1951-1979

\bibitem{WilYu12}
R. Willett and G. Yu, \emph{Higher index theory for certain expanders and Gromov monster groups II},  Adv. in Math. \textbf{229}(2012), 1762-1803


\bibitem{WilYu14}
R. Willett and G. Yu, \emph{Geometric property (T)}, Chin. Ann. Math. \textbf{35}(2014), 761-800

\vspace{2cm}
\textsc{Westf\"alische Wilhelms-Universit\"at M\"unster}

\emph{Address}: Einsteinstrasse 62, 48149 M\"unster

\emph{E-mail address}: jwinkel@uni-muenster.de




 
\end{thebibliography}
\end{document}